\documentclass[a4paper]{article}

\usepackage{amssymb}
\usepackage{amsthm}
\usepackage{amsmath}
\usepackage[cp1250]{inputenc}
\usepackage{graphicx}

\newtheorem{theorem}{Theorem}[section]
\newtheorem{lemma}[theorem]{Lemma}
\newtheorem{rem}[theorem]{Remark}
\newtheorem{definition}[theorem]{Definition}
\theoremstyle{definition}

\DeclareMathOperator{\IM}{Im}

\DeclareMathOperator{\RE}{Re}
\DeclareMathOperator{\spect}{spect}
\DeclareMathOperator{\DOM}{dom}
\DeclareMathOperator{\MID}{mid}

\providecommand{\im}[1]{\IM\left(#1\right)}
\providecommand{\re}[1]{\RE\left(#1\right)}

\providecommand{\jednadruga}{\frac{1}{2}}
\providecommand{\dom}[1]{\DOM\left(#1\right)}

\providecommand{\decay}{order of polynomial decay}

\providecommand{\subspaceH}{H^\prime}

\def\qed{{\hfill{\vrule height5pt width3pt depth0pt}\medskip}}

\begin{document}
\begin{center}
{\Large \bf  Existence of globally attracting  solutions for
one-dimensional viscous Burgers equation with nonautonomous
forcing - a computer assisted proof}

 \vskip 0.5cm
 {\large
Jacek~Cyranka$^{*,\ddag}$ \footnote{Research has been supported by Polish
National Science Centre grant DEC-2011/01/N/ST6/00995.},
Piotr~Zgliczy\'nski$^{*,\dag}$ \footnote{Research has been supported by Polish
National Science Centre grant 2011/03B/ST1/04780}}

 \vskip 0.5cm

{\small$^*$ Institute of Computer Science and Computational Mathematics,
Jagiellonian University}\\
{\small  ul. S. {\L}ojasiewicza 6, 30-348 Krak\'ow, Poland}

 \vskip 0.5cm
 
{\small$^\ddag$ Faculty of Mathematics, Informatics and Mechanics,
University of Warsaw}\\
{\small  Banacha 2, 02-097 Warszawa, Poland}

 \vskip 0.5cm

 {\small$^\dag$ WSB-NLU\\
 ul. Zielona 27, 33--320 Nowy Sacz, Poland}
 
 \vskip 0.5cm

 jacek.cyranka@ii.uj.edu.pl, piotr.zgliczynski@ii.uj.edu.pl

\vskip 0.5cm

\today
\end{center}

\begin{abstract}
  We prove the existence of globally attracting   solutions of the viscous Burgers
  equation with periodic boundary conditions on the interval for some particular choices of
  viscosity and  non-autonomous forcing. The attracting solution is periodic if the forcing is periodic.
    The method is general and
  can be applied to other similar partial differential equations.
  The proof is computer assisted.

\end{abstract}

\section{Introduction}

We present a method of proving the existence of globally
attracting solutions for the viscous Burgers equation with periodic boundary conditions on the interval with a
time-dependent forcing. The attracting solution is periodic, if
the forcing is periodic. It is worth pointing out, that our
method allows to obtain a globally attracting solution in case
of non periodic forcings.
Moreover, the method is general and should be
applicable to other dissipative PDEs with periodic boundary
conditions.

Let us begin with a short review of published results by both of the authors, which 
have been the foundation of our current research. In \cite{C} a method of proving an existence of globally attracting steady-states
for certain class of parabolic PDEs is presented. As an illustration of the method a
detailed case study of the viscous Burgers equation is showed. The method can be summarized in three steps. 
First, construction of global absorbing sets, which are composed from regular functions, and 
such that they absorb any initial condition after a finite time. Second, establishing an existence of
locally attracting steady-state. Third, absorbing sets are showed to 
be mapped into the fixed point's local bassin of attraction by a rigorous numerical
integration procedure. For the purpose of rigorous integration and establishing existence 
of an attracting fixed point we used a topological method of self-consistent bounds, 
developed in the series of articles \cite{ZM, ZAKS, ZNS, Z2, Z3}. 
Now, let us shortly describe the innovation of the presented results. The current results generalize previous works as we deal with 
globally attracting orbits in the nonautonomous case. Establishing new results required from us 
extending the method of self-consistent bounds to the nonautonomous case, and deriving a
new topological principle accompanied by an algorithm of proving existence of locally attracting solution defined on $\mathbb{R}$. 
To prove the attracting solution exists we verify that a time-shift map is contraction in certain space by
construcing an interval enclosure of a so-called isolating segment for discrete semiprocess.
Then we estimate the Lipschitz constant of the time-shift on the calculated enclosure. Observe that 
in \cite{C, ZNS} the forcing was assumed to be constant in time, 
and in \cite{ZM, ZAKS, Z2, Z3} the considered PDEs did not include any external 
forcing at all.

More specifically, in the present paper we present the case study of
the initial value problem with periodic boundary conditions for the Burgers equation 
on the interval
\begin{equation}
    \label{eq:vBEq}
    u_t+u\cdot u_x-\nu u_{xx}=f(t, x),
\end{equation}

 First, as an example result we show
\begin{theorem}
    \label{thm:main1}
    For any $\nu\in[2,2.1]$ and $f \in S_1$, where
\begin{eqnarray*}
    S_1 = \left\{x\mapsto 1.6\cos{2x}-2\sin{3x}+
    \sum_{k=1}^{3}{\beta_k(t)\sin{kx}+\gamma_k(t)\cos{kx}}, \right. \\
\left. \beta_k(t),\ \gamma_k(t)\in\left[-0.03,0.03\right],\forall t\in\mathbb{R} 
\right\},
\end{eqnarray*}
    where $\beta_k(t),\ \gamma_k(t)$ are
    continuous, there exists a classical
    solution (periodic in time when $f$ is periodic in time) of \eqref{eq:vBEq}, 
    defined on $\mathbb{R}$,
    which attracts  any initial data $u_0$ satisfying
      $u_0\in C^4$ and $\int_0^{2\pi}{u_0(x)\,dx}=\pi$. Moreover,  the convergence towards
attracting solution is exponential.
\end{theorem}

Theorem 1.1 in \cite{C} about the existence of globally
attracting fixed points when the forcing is time-independent is a
particular case of Theorem~\ref{thm:main1} with
$\beta_k(t)=\beta_k,\ \gamma_k(t)=\gamma_k$ constant. As we show
in Section~\ref{sec:proof} the computer assisted part of the proof
of this theorem was in fact already accomplished during the proof
of Thm. 1.1 in \cite{C}.

\providecommand{\paperExampleNu}{2}
\providecommand{\paperExampleForcing}{-0.6\sin(x)+0.7\cos(2x)+0.7\sin(2x)-0.8\cos(3x)-0.8\sin(3x)}
\providecommand{\paperExampleAzero}{\pi}
\providecommand{\paperExampleNonautForcing}{\sin(t)\left[-0.6\cos(x)+0.7\cos(2x)+0.7\sin(2x)-0.8\cos(3x)-0.8\sin(3x)\right]}

\providecommand{\paperExampleM}{8}
\providecommand{\paperExampleL}{3.61879e-05}

\providecommand{\paperExampleExecutionTime}{342.34}

\providecommand{\paperExampleEZero}{1.22018}
\providecommand{\paperExampleNu}{2}

\providecommand{\paperExampleFPerturbationDiam}{1e-04}
For the sake of demonstration, we state the next theorem for an explicitly given set of 
nonautonomous forcing functions with a particular dominant part. Our method is not restricted to 
this case, this particular nonautonomous dominant part is a result of setting some parameters in 
our algorithm. Our algorithm is capable to attempt to prove, in principle, any other case in which 
the dominant part is provided by explicit formulas.
\begin{theorem}
    \label{thm:main2}
    For $\nu=\paperExampleNu$ and $f\in S_2$, where
    \begin{multline*}
      S_2 = \left\{x\mapsto\paperExampleForcing+\right.\\
      \paperExampleNonautForcing+\\
      \sum_{k=1}^{3}{\beta_k(t)\sin(kx)+\gamma_k(t)\cos(kx)},\\
      \left.\beta_k(t),\,\gamma_k(t)\in \left[-5\cdot 10^{-5}, 5\cdot 10^{-5}\right],\ \forall t\in\mathbb{R}\right\},
    \end{multline*}
    where $\beta_k(t),\ \gamma_k(t)$ are continuous, there exists a classical
    solution (periodic in time when $f$ is periodic in time) of \eqref{eq:vBEq}, 
    defined on $\mathbb{R}$,
    which attracts  any initial data $u_0$ satisfying $u_0\in C^4$ and
    $\int_0^{2\pi}{u_0(x)\,dx}=\paperExampleAzero$. Moreover,  the convergence towards
attracting solution is exponential.
\end{theorem}

In both theorems we are interested in classical
solutions only. This is the reason, why we do not state the
theorem for more general solutions.

The essential difference between Theorem~\ref{thm:main1} and
Theorem~\ref{thm:main2} is that in Theorem~\ref{thm:main1} the
non-autonomous part of the forcing (the time-dependent part) is a
small perturbation of the autonomous part (the time-independent part).
 Whereas in Theorem~\ref{thm:main2} the norms
of the autonomous, and the non-autonomous part of the forcing are
of the same order of magnitude. Due to this fact the proofs of
both theorems are  based on the slightly different topological
principle. In the proof of Thm.~\ref{thm:main1} we constructed a
trapping isolating segment (a forward invariant set in the extended phases pace),
which is time independent, while in
the proof of Thm.~\ref{thm:main2} we constructed a trapping set (a forward invariant set)
for the time shift along the orbits by $2\pi$ - the period of the
dominant part of the non-autonomous forcing.

Let us comment on the role of the condition
$\int_0^{2\pi}{u_0(x)\,dx}=\pi$ in Theorem~\ref{thm:main1} and
Theorem~\ref{thm:main2}. The condition
$\int_0^{2\pi}{u_0(x)\,dx}\neq 0$, when compared to
$$\int_0^{2\pi}{u_0(x)\,dx}= 0$$ makes the proof significantly
harder numerically, due to the appearance of complex eigenvalues in
the partial derivative of the vector field -- see the numerical
data \eqref{eq:approxEigenv} in Appendix~\ref{sec:numData}.
Therefore for the illustration of our method we
decided to take this more difficult case.

Similar results to ours can be found in literature. In \cite{JKM}
for any $\nu>0$ the authors established the existence of a
globally attracting solution of \eqref{eq:vBEq} periodic in space
and time, under assumption that forcing is periodic in time.
Hence, in this respect, our results for the time-periodic forcing
are significantly less general as we just consider particular
cases of parameters. We believe that even in that case our
approach is of some interest, as we are able to establish the
exponential convergence rate to the attracting solution, while in
\cite{JKM} the authors clearly indicated that they cannot make
such claim and they asked for the convergence rate in one of the
stated problems \cite[Problem 3(i)]{JKM}. The method in \cite{JKM}
appears to be restricted to the scalar equation on one-dimensional
domains, partially due to the use of the maximum principles. In
\cite{Si1} the author established a similar result for
\eqref{eq:vBEq} with a time-periodic forcing proving also the
exponential convergence to the attracting orbit in the periodic
case. The technique used in \cite{Si1} uses heavily the fact that
the Cole-Hopf transformation transforms the Burgers equation to a
linear parabolic equation. This significantly reduces the
applicability of this approach to other PDEs.

The technique we use here is not restricted to some particular
type of equation nor to the dimension one. We need some kind of
'energy' decay as a global property of our dissipative PDEs and
then if the system exhibits an attracting orbit, then we should, in
principle, be able to prove it independently of the system dimensionality. 
Generally speaking, the applicability of our technique to a whole class of
 dPDEs follows from the method of self-consistent bounds properties, 
in \cite{Z3} it is argued how the method of self-consistent bounds 
applies to PDEs with other nonlinearities, in \cite{ZNS} 
the method is applied to the incompressible Navier-Stokes equations with an 
autonomous forcing. In \cite{CZ}, using the averaging principle, we established existence of 
globally attracting orbits, asymptotically for sufficiently large integral of the initial condition,
 for the 1D viscous Burgers and  incompressible 2D Navier-Stokes equations with an 
nonautonomous forcing. 

 At the present state our technique strongly relies
on the existence of good coordinates, the Fourier modes in the
considered example. We hope that the further development of the
rigorous numerics for dissipative PDEs based on other function
bases, e.g.  the finite elements, should allow to treat
also different domains and boundary conditions in the near future.

\subsection{Notation}

Some notation: $\mathbb{R}_{+}=[0,\infty)$, $B(Z,\delta)$ a ball
of size $\delta$ around the set $Z$. $B_n(z,r)$ is a ball in
$\mathbb{R}^n$ with the center $z$ and radius $r$ with the
distance function is known from the context.

We denote by $[x]$ an \textit{interval set}
$[x]\subset\mathbb{R}^n$, $[x]=\Pi_{k=1}^n[x^-_k, x^+_k]$, $[x^-_k, x^+_k]\subset\mathbb{R},\ -\infty<x_k^-\leq x_k^+<\infty$.

For a nonautonomous ODE
\begin{equation}
  x'=f(t,x), \label{eq:notation-Non_auto}
\end{equation}
 where $x \in \mathbb{R}^n$ and $f$ is regular enough to guarantee
uniqueness of the initial value problem $x(t_0)=x_0$ for any $(t_0,x_0)$ for (\ref{eq:notation-Non_auto}),  by
 $\varphi(t_0,t,x)=x(t_0+t)$, where $x(t)$ is a solution (\ref{eq:notation-Non_auto}) with  initial condition  $x(t_0)=x_0$.
 Obviously in each context it will be clearly stated what is the ordinary differential equation generating $\varphi$. We will
 sometimes refer to $\varphi$ as to the local process generated by  (\ref{eq:notation-Non_auto}).


\section{Viscous Burgers equation with periodic boundary conditions on interval}
\label{sec:burgers}

The Burgers equation was proposed in \cite{B} as a mathematical
model of turbulence. There is a significant number of applications
of the Burgers equation, see e.g. \cite{Wh}. We consider the
initial value problem for viscous Burgers equation on the interval
with periodic boundary conditions and \emph{a non-autonomous forcing $F$}, i.e.
\providecommand{\timeInterval}{[t_0, \infty)}
\providecommand{\forcingTimeInterval}{\mathbb{R}}
\begin{subequations}
    \label{eq:burgers}
    \begin{align}
        &u_t(t,x)+u(t,x)\cdot u_x(t,x)-\nu u_{xx}(t,x)=F(t,x),\quad t\in\timeInterval,\ x\in\mathbb{R},\label{eq:burgers1}\\
        &u(t, x)=u(t, x+2\pi),\quad t\in\timeInterval,\ x\in\mathbb{R},\label{eq:burgers3}\\
        &F(t, x)=F(t, x+2\pi),\quad t\in\forcingTimeInterval,\ x\in\mathbb{R},\label{eq:burgers4}\\
        &u(t_0,x)=u_0(x),\quad t_0\in\mathbb{R},\ x\in\mathbb{R},\label{eq:burgers2}
    \end{align}
\end{subequations}
where $\nu >0$.

For the technical purposes we assume that
\begin{equation}
  F(t,x)=f(x) + \widetilde{f}(t,x),
\end{equation}
where $f$ and $\widetilde{f}$ are continuous and
$2\pi$-periodic with respect to $x$ variable, and we define $F(t,x)$ for $t\in\forcingTimeInterval$. Later, we will put
more restrictive conditions on $f$ and $\widetilde{f}$. In fact,
$f$ will be given in an explicit form and for $\widetilde{f}$ we will
demand some bounds.

We will use the Fourier series to study \eqref{eq:burgers}. Let
\begin{equation}
u(t,x)=\sum_{k \in \mathbb{Z}}a_k(t) \exp(ikx).
\end{equation}

It is straightforward to write the problem \eqref{eq:burgers} in
the Fourier basis. We obtain the following infinite ladder
of equations
\begin{equation}
\label{eq:burgers_infinite1}
 \frac{d a_k}{d t}=-i\frac{k}{2}\sum_{k_1\in\mathbb{Z}}{a_{k_1}\cdot a_{k-k_1}}+
        \lambda_k a_k+f_k+\widetilde{f}_k(t),\quad t\in\timeInterval,\ k\in\mathbb{Z},
\end{equation}
where
\begin{subequations}
    \label{eq:burgers_infinite2}
    \begin{align}
        &a_k(t_0)=\frac{1}{2\pi}\int_0^{2\pi}{u_0(x)e^{-ikx}}\,dx,\quad k\in\mathbb{Z},\label{eq:burgers_infinite2}\\
        &f_k=\frac{1}{2\pi}\int_0^{2\pi}{f(x)e^{-ikx}}\,dx,\quad k\in\mathbb{Z},\label{eq:burgers_infinite3}\\
        &\widetilde{f}_k(t)=\frac{1}{2\pi}\int_0^{2\pi}{\widetilde{f}(t,x)e^{-ikx}}\,dx,\quad t\in\forcingTimeInterval,\ k\in\mathbb{Z},\label{eq:burgers_infinite4}\\
        &\lambda_k=-\nu k^2. \label{eq:burgers_infinite5}
    \end{align}
\end{subequations}

 The reality of $u$, $f$ and $\widetilde{f}$ implies
that for $k\in\mathbb{Z}$
\begin{equation}
    a_k=\overline{a_{-k}}, \quad  f_k=\overline{f_{-k}}, \quad
    \widetilde{f}_k(t)=\overline{\widetilde{f}_{-k}(t)}\text{ for }t\in\forcingTimeInterval.
\end{equation}
In view of the above  variables $\{a_k\}_{k\in\mathbb{Z}}$ are not
independent, this motivates the following definition.
\begin{definition}
In the space of sequences $\{a_k\}_{k \in \mathbb{Z}}$,  where
$a_k\in \mathbb{C}$, we will say that the sequence $\{a_k\}$
satisfies the reality condition iff
\begin{equation}
   a_{k}=\overline{a_{-k}}, \quad k \in \mathbb{Z}.  \label{eq:reality-cond}
\end{equation}
We will denote the set of sequences  satisfying
(\ref{eq:reality-cond}) by $R$. It is  easy to see that $R$ is a
vector space over the field $\mathbb{R}$.
\end{definition}

We will assume that the initial condition for
\eqref{eq:burgers_infinite1} satisfies
\begin{equation}
  \label{eq:fixedInt}
  \frac{1}{2\pi}\int_{0}^{2\pi}{u_0(x)\,dx}=\alpha,\quad\text{for a fixed }\alpha\in\mathbb{R}.
\end{equation}
We will require additionally that $f_0=0$, $\widetilde{f}_0(t)=0$ for $t\in\forcingTimeInterval$,
and then \eqref{eq:fixedInt} implies that $a_0(t)$ is constant, namely
\begin{equation}
  \label{eq:fixedA0}
  a_0(t)=\alpha, \quad \forall t \geq t_0.
\end{equation}
\begin{definition}
    \label{def:symmetricGalerkinProjection}
    For any given number $m>0$ {\rm the $m$-th Galerkin projection} of \eqref{eq:burgers_infinite1} is
    \begin{equation}
        \label{eq:symmetricGalerkinProjection}
        \frac{d a_k}{d t}=-i\frac{k}{2}\sum_{\substack{|k-k_1|\leq m\\|k_1|\leq m}}{a_{k_1}\cdot a_{k-k_1}}+\lambda_k a_k+f_k+\widetilde{f}_k(t),\quad t\in\timeInterval,\ |k|\leq m.
    \end{equation}
\end{definition}
Here and further on with a slight abuse of notation we denote
$m$-th Galerkin projection solution's $k$-th mode by $a_k$,
which is the same symbol as the $k$-th mode of the solution of
the full system \eqref{eq:burgers_infinite1}.
Note that the condition \eqref{eq:fixedA0} holds also for all
Galerkin projections \eqref{eq:symmetricGalerkinProjection} as
long as $f_0=0$, and $\widetilde{f}_0(t)=0$ for all
$t\in\forcingTimeInterval$. Also observe that   the reality condition
\eqref{eq:reality-cond} is invariant under all Galerkin
projections \eqref{eq:symmetricGalerkinProjection}, i.e. if
$a_k(t_0)=\overline{a_{-k}}(t_0)$, then $a_k(t)=\overline{a_{-k}}(t)$
for all $t > t_0$ if the solution of
(\ref{eq:symmetricGalerkinProjection}) exists up to that time.

\begin{definition}
    Let ${\nmid}\cdot{\nmid}\colon\mathbb{R}\to\mathbb{R}$ be given by
    \begin{equation*}
      {\nmid}a{\nmid}:=\left\{\begin{array}{ll}|a|&\text{ if }a\neq 0,\\1&\text{ if }a=0.\end{array}\right.
    \end{equation*}
  \end{definition}

\begin{definition}
\label{def:tildeH}
  Let $H$ be the space $l_2(\mathbb{Z},\mathbb{C})$, i.e. $u \in H$ is a sequence $u:\mathbb{Z}\to \mathbb{C}$ such that
  $\sum_{k \in \mathbb{Z}}|u_k|^2 < \infty$  over the coefficient field $\mathbb{R}$.
  The subspace $\widehat{H}\subset H$ is defined by
  \begin{equation*}
    \widehat{H}:=\left\{\{a_k\}\in H\colon\text{there exists }0\leq C<\infty\text{ such that }|a_k|\leq\frac{C}{{\nmid}k{\nmid}^4}\text{ for }k\in\mathbb{Z}\right\}.
  \end{equation*}

  This space is equivalent to the space of sequences having the following weighted $l^\infty$ norm finite
  \begin{equation*}
    \|\{a_k\}\|=\sup_{k\in\mathbb{Z}}{|k|^4|a_k|}
  \end{equation*}
\end{definition}

\begin{definition}
\label{def:H'}
  Let the space $\subspaceH$ be given by
  \begin{equation*}
    \subspaceH:=\widehat{H}\cap R.
  \end{equation*}
\end{definition}

Let us comment on Definitions~\ref{def:tildeH} and \ref{def:H'}. Despite the fact that we are dealing with complex sequences we use as the coefficient field
the set of real numbers, because the reality condition is not compatible with the complex multiplication.

  The choice of the particular subspace $\subspaceH$ is motivated by the fact that the order
  of decay of coefficients $\{a_k\}\in\subspaceH$ is sufficient for the uniform convergence of $\sum{a_ke^{ikx}}$ and every term appearing in
  \eqref{eq:burgers1}. Moreover, in Theorem~\ref{thm:main1} and Theorem~\ref{thm:main2} the attraction property is obtained within the class of $C^4$ functions
  due to the fact that the Fourier expansion of any i.c. $u_0\in C^4$ belongs to $\subspaceH$. For the details see \cite{C}.

\subsection{Absorbing set}

The goal of this section is to establish for the
existence of the forward invariant absorbing set for all Galerkin
projections of (\ref{eq:burgers}), with good compactness
properties. Here, we basically quote the results from \cite{C}
with some  improvements.

\providecommand{\gp}{Galerkin projection of \eqref{eq:burgers_infinite1}}
\providecommand{\gps}{Galerkin projections of \eqref{eq:burgers_infinite1}}

\begin{definition}{\rm\cite[Def. 4.6]{C}}
    \label{def:absorbingSet}
    Let $N_0\geq 0$, $\varphi_n$ be a local process induced by the $n$-th Galerkin projection of
    \eqref{eq:burgers_infinite1}.
    A set $\mathcal{A}\subset \subspaceH$ is called \emph{the absorbing set
    for large \gps}, if for any pair $(t_0,u_0)\in  \mathbb{R} \times \subspaceH$ there exists
    $t_1(u_0)\geq 0$ such that
    for all $n>N_0$ and all $t_0 \in \mathbb{R}$, $t\geq t(u_0)$ holds  $\varphi_n\left(t_0,t, P_nu_0\right)\in P_n
    \mathcal{A}$.
     Moreover, $P_n\mathcal{A}$ is forward invariant for $\varphi_n$.
\end{definition}

Our definition of the absorbing set differs from the standard one of bounded absorbing set,
see for example \cite{FMRT}. There it is
stated for an abstract evolutionary equation and $t_1=t_1(B)$ has to be uniform for any bounded
set $B$, whereas in our case,
we state the definition for the more specific case of (sufficiently)
large Galerkin projections of \eqref{eq:burgers_infinite1} and $t_1$ depends on point $u$,
so we use the notion of point absorbing set. Observe that both mentioned concepts of absorbing sets are
equivalent for a fixed $n$, but as we ask for uniformity in $n>N_0$ we use a weaker concept.
Despite the fact that for the absorbing sets we construct in this work, $t_1$ can be chosen
uniformly for each bounded set $B$, i.e. $t_1=t_1(B)$ we find this stronger requirement unnecessary.

\begin{definition}{\rm\cite[Def. 3.1]{C}}
\label{def:energy}

\emph{Energy} of (\ref{eq:burgers_infinite1}) is given by the formula
\begin{equation}
    \label{eq:energy}
    E(\{a_k\})=\sum_{k\in\mathbb{Z}}{|a_k|^2}.
\end{equation}
Energy of (\ref{eq:burgers_infinite1}) with $a_0$ excluded is given by the formula
\begin{equation}
  \label{eq:energyWithoutZero}
  \mathcal{E}(\{a_k\})=\sum_{k\in\mathbb{Z}\setminus\{0\}}{|a_k|^2}.
\end{equation}

\end{definition}

\providecommand{\Ezero}{\mathcal{E}}
\providecommand{\Fassumption}{
$F_k(t)=\overline{F_{-k}(t)}$, $F_k(t)=0$ for $|k|>J$, and $F_0(t)=0$}
\providecommand{\satisfiesFassumption}{$F_k(t)=f_k+\widetilde{f}_k(t)$
for $t\in\forcingTimeInterval$ satisfies \Fassumption}
\providecommand{\energyBound}{\widetilde{\mathcal{E}}}
\providecommand{\Ef}{\sup_{t\in\forcingTimeInterval}E\left(\left\{F_k(t)\right\}\right)}
\providecommand{\invariantSubspace}{restricted to the invariant
subspace given by $a_k=\overline{a_{-k}}$}
\providecommand{\D}{2^{s-\jednadruga}+\frac{2^{s-1}}{\sqrt{2s-1}}}

The theorem below is a main building block for the
construction of the absorbing set.
\begin{theorem}{\rm Based on \cite[Thm. 3.4]{C}}
  \label{thm:analyticTrappingRegion}
  Assume that \satisfiesFassumption. Let $\{a_k\}_{k\in\mathbb{Z}}\in H$,  $s>0.5$, $E_0=\Ef\nu^{-2} < \infty$,
  $\energyBound>E_0 $, $D=\D$,
  $C>\sqrt{\energyBound}N^s$, $N>\max{\left\{J,\left(\frac{\sqrt{\energyBound}D}{\nu}\right)^{2}\right\}}$.
  Then
    \begin{equation*}
        W(\energyBound, N, C, s)=\left\{ \{a_k\} \in R \ |\ \Ezero(\{a_k\})\leq \energyBound,\ |a_k|\leq\frac{C}{|k|^s}\right\}
    \end{equation*}
    is a \emph{trapping region} (i.e. is forward invariant) for each \gp.
\end{theorem}

The theorem below establishes the existence of a family of
absorbing sets.

\begin{theorem}{\rm Based on \cite[Lemma 4.7]{C}}
    \label{thm:absorbingSet}
     Assume that \satisfiesFassumption.
    Let $\varepsilon>0$,  $E_0=\Ef\nu^{-2} < \infty$,
     $\energyBound>E_0$, $N$ is defined in Thm.~\ref{thm:analyticTrappingRegion}.
    Put
    \begin{eqnarray}
      s_i&=&i/2\text{ for }i\geq 2,\nonumber\\
      D_i&=&2^{s_i-\frac{1}{2}}+\frac{2^{s_i-1}}{\sqrt{2s_i-1}}\text{ for }i\geq 2,\nonumber\\
      C_2&=&\varepsilon+\frac{1}{\nu}\left(\jednadruga\energyBound+\sup_{\substack{0<|k|\leq J\\t\in\forcingTimeInterval}}{\frac{|f_k|+|\widetilde{f}_k(t)|}{|k|}}\right),\label{eq:depOnF1}\\
      C_i&=&\varepsilon+\frac{1}{\nu}\left(C_{i-1}\sqrt{\energyBound}D_{i-1}+\sup_{\substack{0<|k|\leq J\\t\in\forcingTimeInterval}}{|k|^{s_i-2}\left(|f_k|+|\widetilde{f}_k(t)|\right)}\right)\text{ for }i>2,\label{eq:depOnF2}.
    \end{eqnarray}\\
    Then for all $i\geq 2$, and $\widetilde{C}_i>\sqrt{\energyBound}N^{s_i}$
    \begin{equation*}
        H\supset \mathcal{W}_i\bigl(\energyBound,C_i,\varepsilon\bigr):=\left\{\{a_k\}_{k\in\mathbb{Z}}\in R\ |\ \Ezero(\{a_k\}_{k\in\mathbb{Z}})\leq\energyBound,\ |a_k|\leq\frac{C_i}{|k|^{s_i}}\right\}\bigcap W(\energyBound,N,\widetilde{C}_i, s_i),
    \end{equation*}
    is an \emph{absorbing set} for large \gps.
\end{theorem}
The absorbing sets obtained in the above  theorem, contrary to
\cite[Lemma 4.7]{C}, does not depend on $\alpha$
\eqref{eq:fixedInt}.  As an
consequence of Theorem~\ref{thm:absorbingSet}, and some
improvements of the algorithms presented in \cite{C}, we are not
anymore constrained with large $\alpha$ values. We managed to
prove some example theorems for cases with large $\alpha$ values,
and the results are presented in Table~\ref{table2}.

The intersection with the trapping isolating
segment $W$ is required to ensure the obtained set is forward
invariant in time. The proof of Theorem~\ref{thm:absorbingSet}
follows the scheme of the proof of \cite[Lemma 4.7]{C}, however
the following auxiliary lemma is required.
Precisely, for the sake of proving Theorem~\ref{thm:absorbingSet},
\cite[Lemma 4.4]{C} should be replaced by Lemma~\ref{lem:akbk2} below.

\begin{lemma}
    \label{lem:akbk2}
    Assume that \satisfiesFassumption.
    Let  $\subspaceH\supset W$ be trapping region (i.e. is forward invariant) for all
    \gps, such that $P_n W \subset W$ for all $n$.

    Assume that $C_a\geq 0$, $s_a>0.5$ are numbers such that
    \begin{equation}
      \label{assum1}
      |a_k|\leq\frac{C_a}{|k|^{s_a}} \quad \text{ for }|k|>0\text{, and for all }a\in W.
    \end{equation}
    Assume that $C_\mathcal{N}\geq 0$, $s_\mathcal{N}\geq s_a-1$ are numbers such that
    \begin{equation*}
        \left|\mathcal{N}_k(a)\right|\leq\frac{C_\mathcal{N}}{|k|^{s_\mathcal{N}}}
        \quad \text{ for }|k|>0, \text{and for all }a\in W,
    \end{equation*}
    where
    \begin{equation}
      \label{assum2}
      \mathcal{N}_k(a)=-i\frac{k}{2}\sum_{k_1\in\mathbb{Z}\setminus\{0,k\}}{a_{k_1}\cdot a_{k-k_1}}.
    \end{equation}

    Then for any $\varepsilon>0$ there exists a finite time $\hat{t}\geq 0$ such that for all $l>0$ and $t\geq\hat{t}$,
    any $a(t_0+t)$ -- the solution of $l$-th Galerkin projection of \eqref{eq:burgers_infinite1} such that $a(t_0)\in P_l(W)$, satisfies
    \begin{equation}
            |a_k(t_0+t)|\leq\frac{C_b+\varepsilon}{|k|^{s_b}}\text{ for }0<|k|\leq l,  \label{eq:akt-exp-small}
    \end{equation}
    where
    \begin{equation}
      C_b=\left(C_\mathcal{N}+\sup_{\substack{0<|k|\leq J\\ t \in \mathbb{R}}}{\left\{|F_k(t)||k|^{s_\mathcal{N}}\right\}}\right)/\nu, \quad
    s_b=s_\mathcal{N}+2.  \label{eq:cbsb}
    \end{equation}
\end{lemma}
\paragraph{Proof}

 Let us fix the Galerkin projection
dimension $l>0$.

We consider the initial value problem for the l-th Galerkin projection  of \eqref{eq:burgers_infinite1} with the initial condition $a(t_0)=a_0\in P_lW \cap \subspaceH$.

 Using the reality condition \eqref{eq:reality-cond} we obtain
\begin{multline*}
  \frac{d |a_k|^2}{d t}=\frac{d a_k}{d t}\cdot a_{-k}+a_k\cdot \frac{d a_{-k}}{d t}=\\
  \left((\lambda_k-ika_0)a_k+\mathcal{N}_k(a)+F_k(t)\right)a_{-k}+\left((\lambda_k+ika_0)a_{-k}+\mathcal{N}_{-k}(a)+F_{-k}(t)\right)a_{k}=\\
  2\lambda_k|a_k|^2+\left(\mathcal{N}_k(a)+F_k(t)\right)a_{-k}+\left(\mathcal{N}_{-k}(a)+F_{-k}(t)\right)a_{k}.
\end{multline*}

 From the reality condition for $a$, $\mathcal{N}$ and $F$ we obtain
\begin{equation*}
  \frac{d |a_k|^2}{d t}\leq 2\lambda_k|a_k|^2+2\left(\sup_{t\in\mathbb{R}, u \in W}\left|\mathcal{N}_{-k}(u)+
F_{-k}(t)\right|\right) \left|a_k\right|,
\end{equation*}
hence
\begin{equation}
  \frac{d |a_k|}{d t}\leq \lambda_k|a_k|+\sup_{t\in\mathbb{R}, u \in W}\left|\mathcal{N}_{-k}(u)+
F_{-k}(t)\right|\quad\text{for}\ |a_k|>0.\label{eq:normEstimateN}
\end{equation}
Let
\begin{equation}
  b_k=\sup_{t \in \mathbb{R}, u \in W}\left|\mathcal{N}_{-k}(u)+F_{-k}(t)\right|/(-\lambda_k).  \label{eq:bk}
\end{equation}
From \eqref{assum1}, \eqref{assum2}, \eqref{eq:cbsb} it follows that
\begin{equation}
  b_k \leq \frac{C_b}{|k|^{s_b}},\quad |k|>0. \label{eq:b-estm}
\end{equation}
From \eqref{eq:normEstimateN} it follows that for $t \geq 0$ holds
\begin{equation*}
  \left|a_k(t_0+t)\right|\leq\left(\left|a_k(t_0)\right|-b_k\right)e^{\lambda_{k}t}+b_k,\text{ for }|k|>0.
\end{equation*}
From \eqref{eq:b-estm}, and \eqref{assum1} we obtain for $t \geq 0$ and $|k|>0$
\begin{equation*}
  \left|a_k(t_0+t)\right|\leq\left(\frac{C_a}{|k|^{s_a}}-\frac{C_b}{|k|^{s_b}}\right)e^{\lambda_{k}t}+\frac{C_b}{|k|^{s_b}},\text{ for }|k|>0,
\end{equation*}
 We would like to find $\hat{t}$ such that for $t \geq \hat{t}$ condition \eqref{eq:akt-exp-small} is satisfied. It is easy to see that
 this is implied by the following inequality, which should be satisfied for $|k| \geq 1$
 \begin{equation}
   C_a|k|^{s_b-s_a}e^{\lambda_{k}t} \leq \epsilon.
 \end{equation}
 Observe that $s_b > s_a$. Let us fix $n \in\mathbb{Z}_+$ such that $n > (s_b-s_2)/2$. We have for $t>0$ and any $|k| \geq 1$
 \begin{eqnarray*}
   C_a|k|^{s_b-s_a}e^{\lambda_{k}t} = \frac{C_a |k|^{s_b-s_a}}{e^{\nu |k|^2 t}} \leq \frac{C_a |k|^{s_b-s_a}}{(\frac{(\nu |k|^2 t)^n}{n!})} \leq \frac{n! C_a}{\nu^n t^n} < \epsilon
 \end{eqnarray*}
 for $t\geq \hat{t}$, $\hat{t}$ is large enough (independent of the dimension of the Galerkin projection, but depending on the set $W$). This finishes the proof of condition \eqref{eq:akt-exp-small}.

\qed


\section{Topological theorems}
\label{sec:topthm}

In this section we state two topological theorems, which are used
to obtain the attracting orbits. It is based on forward invariant sets (trapping regions) and the Brouwer theorem.
We will use the terminology of  \emph{the isolating segment} introduced by R.
Srzednicki (see \cite{S1,SW}) and local processes.

\subsection{Semiprocesses and nonautonomous differential equations}
\label{subsec:semi-proc}

We start with introducing the notion of a {\em local semiprocess}
which formalizes the notion of a continuous family of local
forward trajectories in an extended phase--space.

\begin{definition}
\label{def:semi} Assume that $X$ is a topological space and
$\varphi :D\rightarrow X$ is a continuous mapping, $D \subset
\mathbb{R} \times \mathbb{R}_{+}\times X$ is an open set. We will
denote by $\varphi_{(\sigma,t)}$ the function
$\varphi(\sigma,t,\cdot)$.

$\varphi$ is called a {\em local semiprocess} if the following
conditions are satisfied
\begin{description}
\item[(S1)]
$\forall \sigma \in \mathbb{R}$, $x\in X$ : $\{ t\in
\mathbb{R}_{+}:(\sigma,x,t)\in D \}$ is an interval,
\item[(S2)]
$\forall \sigma \in \mathbb{R} : \varphi_{(\sigma,0)}={\rm
id}_{X}$
\item[(S3)]
$\forall \sigma \in \mathbb{R}, \forall s,t \in \mathbb{R}_+ :
\varphi_{(\sigma,s+t)}=\varphi_{(\sigma+s,t)}\circ
\varphi_{(\sigma,s)}$,
\end{description}
If $D=\mathbb{R} \times \mathbb{R}_+ \times X$, we call $\varphi$
a {\em (global) semiprocess}. If $T$ is a positive number such
that
\begin{description}
\item[(S4)]
$\forall \sigma, t \in
\mathbb{R}_+:\varphi_{(\sigma+T,t)}=\varphi_{(\sigma,t)}$
\end{description}
we call $\varphi$ a {\em $T$-periodic local semiprocess}.
\end{definition}

A local semiprocess $\varphi$ on $X$ determines a local semiflow
$\Phi$ on $\mathbb{R}\times X$ by the formula
\begin{equation}
 \Phi_{t}(\sigma,x)=(\sigma+t,\varphi_{(\sigma,t)}(x)).  \label{eq:semiflow}
\end{equation}
In the sequel we will often call  the first coordinate in the
extended phase space $\mathbb{R} \times X$  {\em a time}.

Let $\varphi$ be  a  local semiprocess  and let $\Phi$ be a local
semiflow associated to $\varphi$. It follows by $(S1)$ and $(S2)$
that for every $z=(\sigma,x)\in \mathbb{R}\times X$ there is an
$0<\omega_{z}\leq +\infty$ such that $(\sigma,t,x)\in D$ if and
only if $0\leq t <\omega_{z}$. Let $x\in X$, $\sigma \in
\mathbb{R}$, then {\em a left solution} through $z=(\sigma,x)$ is
a continuous map $v: (a,0]\rightarrow \mathbb{R}\times X$ for some
$a\in [-\infty,0)$ such that:
\begin{description}
\item[(I)]
$v(0)=z$,
\item[(II)]
for all $t\in (a,0]$ and $s>0$ with $s+t\leq 0$ it follows that
$s<\omega_{v(t)}$ and $\Phi_{s}(v(t))=v(t+s)$.
\end{description}
If $a=-\infty$ then we call $v$ {\em a full left solution}. We can
extend a left solution through $z$ onto $(a,0]\cup[0,\omega_{z})$
by setting $v(t)=\Phi_{t}((\sigma,x))$ for $0\leq t<\omega_{z}$,
to obtain {\em a solution through $z$}. If $a=-\infty$ and
$\omega_{z}=+\infty$, $v$ is called {\em a full solution}. If for
each $x \in X$ $\omega_x=\infty$, then we will say that $\Phi$ is
a \emph{global semiprocess}.

\begin{rem}
The differential equation
\begin{equation}
   \dot{x}=f(t,x)  \label{eq:*}
\end{equation}
such that $f$ is regular enough to guarantee the uniqueness for
the solutions of the Cauchy problems associated to (\ref{eq:*})
generates a  local process as follows: for $x(t_{0},x_{0};\cdot)$
the solution of (\ref{eq:*}) such that
$x(t_{0},x_{0};t_{0})=x_{0}$ we put
\begin{equation}
\varphi_{(t_{0},\tau)}(x_{0})=x(t_{0},x_{0};t_{0}+\tau).
\end{equation}
 If $f$
is $T$-periodic with respect to $t$ then $\varphi$ is a
$T$-periodic local process. In order to determine all $T$-periodic
solutions of equation (\ref{eq:*}) it suffices to look for fixed
points of $\varphi_{(0,T)}$.
\end{rem}

\subsection{Trapping isolating segments}
\label{subsec:trapp-iso-seg}
 We use the following notation: by
$\pi_{1}:\mathbb{R}\times X \rightarrow \mathbb{R}$ and
$\pi_{2}:\mathbb{R} \times X \rightarrow X$ we denote the
projections and for a subset $Z\subset \mathbb{R} \times X$ and
$t\in \mathbb{R}$ we put
$$
Z_{t}= \{ x \in X: (t,x)\in Z \}.
$$

Now we are going to state the definition of the trapping isolating
segment, which is a modification of the notions of $T$-periodic
isolating segment and periodic isolating segment over $[0,T]$
from \cite{S1,SW}.

\begin{definition}
\label{def:setper} We will say that a set $Z \subset {\mathbb{R}}
\times X$ is $T$--periodic, iff $Z_{T+t}=Z_t$ for every $t \in {\mathbb{R}}$.
\end{definition}

We remark that according to the following definition $W(\energyBound,N, C, s)$ from
Theorem~\ref{thm:analyticTrappingRegion} is a \emph{trapping isolating segment}.

\begin{definition}
\label{def:blok} Let $W \subset {\mathbb{R}}\times X$. We call $W$
a {\em trapping isolating segment} for the  global semiprocess
$\varphi$ if:

\begin{description}
\item[(i)] $W \cap ([t_1,t_2]\times X)$  is a  compact set for any $t_1,t_2 \in \mathbb{R}$

\item[(ii)]
for every $\sigma \in \mathbb{R}$, $x\in \partial W_{\sigma}$
there exists $\delta>0$ such that for all $t\in (0,\delta)$
$\varphi_{(\sigma,t)}(x)\in {\rm int} W_{\sigma+t}$.

\end{description}

\end{definition}

Further we will need a notion of the trapping isolating segment
for a differential inclusion
\begin{equation}
   x'(t) \in f(t,x(t)) + [\delta], \label{eq:ode-incl}
\end{equation}
where $x(t) \in \mathbb{R}^n$ and $[\delta] \subset \mathbb{R}^n$.
\begin{definition}
Let $f$ be $C^1$ with respect to $x$, $\frac{\partial f}{\partial x}$
and $f$ be continuous with respect to $(t,x)$. We will say that $W
\subset {\mathbb{R}}\times \mathbb{R}^n$ is a trapping isolating
segment for (\ref{eq:ode-incl}) iff for any function $g:\mathbb{R}
\times \mathbb{R}^n \to \mathbb{R}^n$,  $C^1$ with respect to $x$,
$\frac{\partial g}{\partial x}$ and $g$ continuous with respect to
$(t,x)$, such that $g(t,x) \in [\delta]$, the set $Z$ is a trapping
isolating segment for the semiprocess induced by
\begin{equation}
   x'(t) = f(t,x(t)) + g(t,x(t)). \label{eq:ode-sel-incl}
\end{equation}

\end{definition}

\begin{theorem}
\label{thm:traping-per}
 Assume that $W$ is a T-periodic trapping
isolating segment for $T$-periodic global semiprocess $\varphi$
and $W_0$ is homeomorphic to $\overline{B}_n(0,1)$.

Then $W$ contains a $T$-periodic orbit.
\end{theorem}
\textbf{Proof:} Let $P$ be the map given by the time
shift by $T$. $P$ is defined on $W_0$ and we have $P(W_0) \subset
W_0$. The Brouwer theorem implies the existence of $x \in W_0$,
such that $P(x)=x$, which give rise to a $T$-periodic orbit.
 \qed

\begin{theorem}
\label{thm:trap-exists-orbit}
 Assume $W$ is a  trapping isolating segment for
a global semiprocess $\varphi$ induced by a non-autonomous ODE
\begin{equation}
   x'=f(t,x), \quad f \in C^1(\mathbb{R}\times
   \mathbb{R}^n,\mathbb{R}^n).
\end{equation}

Then there exists $x \in W_0$, such that there exists a full
orbit (forward and backward) through $x$  contained in $W$.
\end{theorem}
\textbf{Proof:}
 Each forward orbit starting from $W_0$ is
contained in $W$. Therefore it is enough to prove the existence of
full backward orbit in $W$.

It is easy to see that for any $l \in \mathbb{N}$ there exists
$v_l:[-l,0] \to \mathbb{R}^n$ an orbit of our semiprocess
contained in $W$.

We would like to show that we can chose a subsequence
$\{u_{l_k}\}$ such that $u_{n_l}$ is converging locally uniformly
 on $(-\infty,0]$ to some full backward orbit $v$.

Let us fix any $k \in \mathbb{N}$. Observe that for $l>k$  $v_l$
is defined on $[-k,0]$ and are contained in $W \cap ([-l,0] \times
\mathbb{R}^n)$, which is a compact set. Therefore there exists
$M>0$ such that
\begin{equation}
  |f(t,x)| \leq M, \qquad (t,x) \in W \cap ([-k,0] \times
\mathbb{R}^n).
\end{equation}
Therefore
\begin{equation}
  |v_l'(t)| \leq M , \quad t \in [-k,0].
\end{equation}
This shows that functions $\{v_l:[-k,0] \to  \mathbb{R}^n\}$ are
equicontinuous  and contained in a bounded set $\pi_x (W \cap
([-l,0] \times \mathbb{R}^n)))$. It follows from the Ascoli-Arzela
Theorem that we can chose subsequence $\{v_{l_m}\}$ which is
uniformly converging on $[-k,0]$ to some continuous function
$\overline{v}_k:[-k,0] \to \mathbb{R}^n$, which is an orbit of the
semiprocess.

Now let us consider the following procedure: assume that we have a
subsequence of solutions $v_{l_i}$ converging uniformly to
$\overline{v}_k:[-k,0] \to \mathbb{R}^n$. From that sequence we
can chose a subsequence which will be uniformly converging to
$\overline{v}_{2k}:[-2k,0] \to \mathbb{R}^n$ and then we find a
subsequence converging on $[-2^2k,0]$ to $\overline{v}_{2^2k}$ and
so on. From all these nested subsequences  $\{u_{l_i}\}$
 by
choosing diagonal elements $\{u_{l_k}\}$ we obtain a sequence,
which is converging uniformly on each compact interval $[-k,0]$ to
$\overline{v}$, such that $\overline{v}(t)=\overline{v}_k(t)$ for
$t \in [-k,0]$. From the continuity of $\varphi$ it follows easily
that $\overline{v}$ is a full backward orbit of $\varphi$.
Obviously, $\overline{v}$ is contained in $W$.
\qed

\subsection{Discrete semiprocesses -  iterations of maps}
\begin{definition}
Assume that we have an indexed family of continuous maps
$\{f_i:\mathbb{R}^n \to \mathbb{R}^n \}_{i \in \mathbb{Z}}$. We
define a map $\varphi: \mathbb{Z}\times \mathbb{Z}_+ \times
\mathbb{R}^n \to \mathbb{R}^n$ by
\begin{equation}
 \varphi(i_0,i,x)=
  \begin{cases}
    x & \text{if $i=0$}, \\
   f_{i_0+i-1} \circ \cdots f_{i_0+1} \circ f_{i_0}(x) & \text{otherwise}.
  \end{cases}
\end{equation}
$\varphi$  we will be called \emph{a discrete semiprocess}.

For $T \in \mathbb{Z}_+$ we say that $\varphi$ is $T$-periodic, if
$f_{i+T}=f_i$ for all $i \in \mathbb{Z}$.
\end{definition}

Analogously with the continuous case define the notion of the
forward and backward orbit for a discrete semiprocess.

\begin{definition}
Consider a set $W=\Pi_{k \in \mathbb{Z}} W_k$. It will be called
\emph{a trapping isolating segment} for the discrete semiprocess
$\varphi$ if the following conditions are satisfied
\begin{description}
\item[(i)] $W_k$  is a  compact set for any $k \in \mathbb{Z}$

\item[(ii)]
for every $k \in \mathbb{Z}$
\begin{equation}
  \varphi(k,1,W_k) \subset {\rm int} W_{k+1}.
\end{equation}
\end{description}
For $T \in \mathbb{Z}_+$ we say that $W$ is $T$-periodic if
$W_k=W_{T+k}$ for all $k \in \mathbb{Z}$.
\end{definition}

We now establish  discrete versions of theorems from
Section~\ref{subsec:trapp-iso-seg}.

\begin{theorem}
\label{thm:discrete-traping-per}
 Assume that $W$ is a T-periodic trapping
isolating segment for a discrete $T$-periodic  semiprocess
$\varphi$ and $W_0$ is homeomorphic to $\overline{B}_n(0,1)$.

Then $W$ contains a $T$-periodic orbit.
\end{theorem}
The proof is the same as in the continuous case and will be
omitted.

\begin{theorem}
\label{thm:discrete-trap-exists-orbit}
 Assume $W$ is a  trapping isolating segment for
a  discrete semiprocess $\varphi$.

Then there exist $x \in W_0$  and a full orbit (forward and
backward) through $x$  contained in $W$.
\end{theorem}
The proof of this theorem uses the same idea as the proof of
Theorem~\ref{thm:trap-exists-orbit}, but in the discrete case
there is no need for the equicontinuity and the Ascoli-Arzela
theorem.


\section{The bounds for the Lipschitz constant for the time evolution of dissipative PDEs}
\label{sec:Lip-const-flow}

\subsection{Basic theorem on logarithmic norms and ODEs}

Consider now the differential equation
\begin{equation}
  x'=f(t,x), \label{eq:odelogn}
\end{equation}
where $f$ and $\frac{\partial f}{\partial x}$ are  is continuous.

 By $\|x\|$ we denote a fixed arbitrary norm in $\mathbb{R}^n$. Let $\mu: \mathbb{R}^{n \times n} \to \mathbb{R}$
be the \emph{logarithmic norm} of $A$ induced by norm $\|\cdot\|$, which  was introduced independently by Dahlquist \cite{D} and Lozinskii \cite{L}  (see also \cite{HNW,KZ} and references given there)
 \begin{equation}
   \mu(A) = \lim_{h \to 0^+} \frac{\|1 + h A\| - 1}{h}.  \label{eq:def-logn}
 \end{equation}
 Observe that $\mu(A)$ is not a norm, as it can be negative. 
 
 It was introduced, because  gives us the bound for the Lipschitz constant of the the time shift
 by $h >0$  of the flow for trajectories contained in a convex compact set $W$ for $t \in [t_0,t_0+h]$ in the form
 \begin{equation} 
 L=\exp\left(h \max_{(t,x) \in [t_0,t_0+h] \times W} \mu\left(\frac{\partial f}{\partial x}(t,x)\right)\right).
 \end{equation}
 Observe that might 
 to be less than one for attracting orbits (because the logarithmic norm can be negative). This should be contrasted with a more standard bound coming from the Gronwall inequality in the form
 \begin{equation}
  L=\exp\left(h \max_{(t,x) \in [t_0,t_0+h] \times W} \left\|\frac{\partial f}{\partial x}(t,x) \right\|\right),
 \end{equation}
   which never can be less than one.
   
A good illustration of the above phenomenon is a linear 1D equation 
\begin{equation*}
x'=-x=f(t,x).
\end{equation*}  
In this case we have the  norm $\|x\|=|x|$. Since $\frac{\partial f}{\partial x}=[-1]$, so $\mu\left(\frac{\partial f}{\partial x} \right) = -1 <0 $ and we obtain a correct bound for the Lipschitz constant, $L$, of the time shift by $h >0$ given by
\begin{equation}
  L=e^{-h}.
\end{equation}
 
Depending on the norm $\|\cdot\|$ the formula for $\mu(A)$ differ. Let us list it for several popular norms
\begin{itemize}
 \item for euclidian norm, 
    \begin{equation}  
      \mu(A) = \max_{\lambda} \left\{ \lambda \in \sigma \left( \frac{A+A^T}{2}   \right) \right\}, \label{eq:eucl-ln}
    \end{equation}
    where $\sigma(M)$ is the spectrum of the matrix $M \in \mathbb{R}^{n \times n}$,
 \item for $\|x\|_\infty:=\max_{i=1,\dots,n} |x_i|$,
 \begin{equation}
   \mu(A)=a_{ii} - \sum_{j,j\neq i} |a_{ij}|, \label{eq:infty-ln}
 \end{equation}  
  \item for $\|x\|_1:=\sum_i |x_i|$,
 \begin{equation}
   \mu(A)=a_{ii} - \sum_{j,j\neq i} |a_{ji}|. \label{eq:l1-ln}
 \end{equation}  
\end{itemize}
In our work with attracting orbits we will always try to change coordinates so that the diagonal of $\frac{\partial f}{\partial x}$ dominates and has only negative entries. In this way
we obtain that the logarithmic norm is negative (see formulas (\ref{eq:eucl-ln},\ref{eq:infty-ln},\ref{eq:l1-ln}) ).

The following theorem is a precise statement on how to obtain the Lipschitz constant for the flow using the logarithmic norm.
It was proved in \cite[Th. I.10.6]{HNW} (we use
a different notation ).
\begin{theorem}
\label{thm:lognorm} Let $y:[t_0,t_0+T] \to \mathbb{R}^n$ be a continuous  piecewise
$C^1$ function and $x:[t_0,t_0+T]  \to \mathbb{R}^n $ be a solution of (\ref{eq:odelogn}).

 Suppose that  the following estimates hold:
\begin{eqnarray*}
  \mu\left(\frac{\partial f}{\partial x}(t,\eta)\right) &\leq& l(t),\quad \mbox{ for $\eta \in [y(t),x(t)] $} \\
  \left\| \partial_+ y(t) - f(t,y(t)) \right\| &\leq& \delta(t).
\end{eqnarray*}
where $\partial_+ y(t)=\lim_{h \to 0^+} \frac{y(t+h)-y(t)}{h}$, i.e. it is the right derivative of $y$ at $t$.

Then for $t > t_0$ holds
\begin{displaymath}
 \| x(t) - y(t)  \| \leq e^{L(t)}\left( \|y(t_0) - x(t_0) \| + \int_{t_0}^{t}e^{-L(s)}\delta(s)ds  \right),
\end{displaymath}
where $L(t)=\int_{t_0}^{t}l(s)ds$.
\end{theorem}

\subsection{Lipschitz constants for the time evolution}
\label{subsec:LipConstFlow}
 We consider a nonautonomous problem
\begin{equation}
  \frac{d a}{dt}=G(t,a), \label{eq:nonAuto}
\end{equation}
where $a \in \mathbb{R}^d$, $G:\mathbb{R} \times \mathbb{R}^d \to
\mathbb{R}^d$ is $C^1$ with respect to $a$ and $G,\frac{\partial
G}{\partial a}$ are continuous.

Let $\varphi(t_0,t,x)$ be a local process induced by (\ref{eq:nonAuto}).

From Theorem~\ref{thm:lognorm} we can easily obtain
the following lemma, which   expresses the Lipschitz constant for the
semiprocess induced by (\ref{eq:nonAuto}) in terms of logarithmic
norms of $DG$ along the trajectory.
\begin{lemma}
\label{lem:lip1} Let  $t_0 < t_1 < t_2 < \dots < t_n$, $[x_i]
\subset \mathbb{R}^d$ for $i=0,\dots,n$,
 $[W_i] \subset \mathbb{R}^d $ for $i=1,\dots,n$ be convex sets and
 $l_i \in \mathbb{R}$ are
 such that
 \begin{eqnarray*}
   \varphi(t_{i-1},[0,t_i - t_{i-1}],[x_{i-1}]) &\subset& [W_{i}],
   \quad i=1,\dots,n \\
   \sup_{(t,a) \in  [t_{i-1},t_i] \times [W_i]} \mu\left(\frac{\partial G}{\partial a}(t,a)\right) &\leq& l_i, \quad i=1,\dots,n \\
   \varphi(t_{i-1},t_i - t_{i-1},[x_{i-1}]) &\subset& [x_{i}],
   \quad i=1,\dots,n
 \end{eqnarray*}
 Then for any $z_1,z_2  \in [x_0]$  holds
 \begin{equation}
  \| \varphi(t_0,t_n-t_0,z_1) - \varphi(t_0,t_n-t_0,z_2) \| \leq \exp\left(\sum_{i=1}^n
  l_i(t_{i}-t_{i-1})\right)\|z_1-z_2\|.
 \end{equation}
\end{lemma}

 In the context
of the above lemma we need to allow for the changes of norms.
We will assume that for $t \in [t_{i-1},t_i]$ we have a norm
$\|\cdot\|_i$. We also assume that there exists norm $\|\cdot\|_0$
just for $t=t_0$. Therefore for $t_i$, $i=1,\dots,n$ we have two
norms. We assume that
\begin{equation}
  \|x\|_i \leq P_{i\mapsto i+1} \|x\|_{i+1}, \quad i=0,\dots, n-1.
\end{equation}

In that context we reformulate the above lemma as follows
\begin{lemma}
\label{lem:lip2} Let $t_0 < t_1 < t_2 < \dots < t_n$,  $[x_i]
\subset \mathbb{R}^d$ for $i=0,\dots,n$,
  $[W_i] \subset \mathbb{R}^d$ for $i=1,\dots,n$ be convex sets,
 $l_i \in \mathbb{R}$ are such that
 \begin{eqnarray*}
   \varphi(t_{i-1},[0,t_i - t_{i-1}],[x_{i-1}]) &\subset& [W_{i}],
   \quad i=1,\dots,n \\
   \sup_{(t,a) \in [t_{i-1},t_i] \times [W_i]} \mu_i\left(\frac{\partial G}{\partial a}(t,a)\right) &\leq& l_i, \quad i=1,\dots,n \\
   \varphi(t_{i-1},t_i - t_{i-1},[x_{i-1}]) &\subset& [x_{i}],
   \quad i=1,\dots,n
 \end{eqnarray*}
 Then for any $z_1,z_2\in [x_0]$ holds
 \begin{equation*}
  \| \varphi(t_0,t_n-t_0,z_1) - \varphi(t_0,t_n-t_0,z_2) \|_n \leq  L \|z_1-z_2\|_0.
 \end{equation*}
 where
 \begin{equation}
   L= \Pi_{i=1}^n \left( \exp(l_i (t_i-t_{i-1})) P_{i-1 \mapsto i} \right)
 \end{equation}
\end{lemma}


\section{Tools for attracting orbits}
\label{sec:top-lip-const}

In this section we consider (\ref{eq:nonAuto}) and we assume that
$G$ satisfies the regularity assumptions from Section
\ref{sec:Lip-const-flow}, i.e. $G$ and $\frac{\partial G}{\partial
a }$ are continuous.

\begin{theorem}
\label{thm:trap--attracting-exists-orbit}
 Assume $W$ is a convex trapping isolating segment for
a global semiprocess $\varphi$ induced by (\ref{eq:nonAuto}).

Assume that
\begin{equation}
\sup_{(t,z) \in W} \mu\left(\frac{\partial G}{\partial
z}(t,z)\right) \leq l.
\end{equation}

Then there  exists a full orbit $v$ for $\varphi$, such that for
any $(t_0,z_0)$ in $W$ and $t >0$
\begin{equation}
 \|\varphi(t_0,t,z_0) - v(t_0+t)\| \leq \exp (l t) \|v(t_0) -
  z_0\|.  \label{eq:lip-const-attr}
\end{equation}
If $l<0$, then the orbit $v$ attracts all other points in $W$.

If  $W$ is a T-periodic trapping isolating segment with $W_0$
homeomorphic to $\overline{B}_n(0,1)$ and $\varphi$ is
$T$-periodic global semiprocess, then the orbit $v$ is
$T$-periodic.
\end{theorem}
\textbf{Proof:}
The existence of the full orbit contained in $W$ follows immediately from Theorem~\ref{thm:trap-exists-orbit}.  In the case of $T$-periodic
semiprocess and trapping isolating segment the existence of the periodic orbit follows from Theorem~\ref{thm:traping-per}.

To obtain  (\ref{eq:lip-const-attr}) observe that $v(t_0+t)=\varphi(t_0,t,v(t_0))$ and we use Theorem~\ref{thm:lognorm} with $x(t)=\varphi(t_0,t,z_0)$
and $y(t)=v(t_0+t)$ for arbitrary $t>0$. Observe that in this situation $\delta(t)=0$, because $y(t)$ is a solution of(\ref{eq:nonAuto}).
\qed

The  theorem given above will be used in the context of the
time-independent isolating segment. The next  theorem we want to
apply in the situation, when finding of an isolating segment for
which the logarithmic norm is negative  appears to be  very
difficult, but it turns out the time shift by the period of the
dominant non-autonomous part has a ball which is mapped into
itself.

Let us fix $T>0$. We define the discrete semiprocess by setting
\begin{equation}
  g_i(x)=\varphi(iT,T,x), \label{eq:T-time-shift}
\end{equation}
i.e. this a time shift by $t=T$ from the section $t=iT$ to
$t=(i+1)T$.

\begin{theorem}
\label{thm:discrete-trap-attr-exists-orbit}
 Assume $W$ is a  trapping isolating segment for
discrete semiprocess (\ref{eq:T-time-shift}).

Assume that there exists  compact and convex set $Z \subset
\mathbb{R}^d$ and $L,B \in \mathbb{R}$ such that for $i\in
\mathbb{Z}$ holds
\begin{eqnarray}
  \sup_{x \in W_i} \|Dg_i(x)\| &\leq& L , \\
   \varphi(iT,[0,T], W_i) &\subset& Z, \\
   \sup_{(t,z) \in \mathbb{R} \times
  Z} \mu\left(\frac{\partial G}{\partial
z}(t,z)\right) &\leq& B. \label{eq:muG-estm2}
\end{eqnarray}

Then there  exists a full orbit $v$ for $\varphi$,
$l=\frac{\ln L}{T}$ and $C=\max(1,\exp(B T)) \cdot \max(1,\exp(-lT))$, such that for any
$(k,z)$ in $W$ and $t >0$
\begin{equation}
 \|\varphi(kT,t,z) - v(kT+t)\| \leq C\exp (l t) \|v(kT) -
  z\|.  \label{eq:lip-iter}
\end{equation}
If $l<0$, then the orbit $v$ attracts all other points in $W$.

If   $W$ is k-periodic  for some $k \in
\mathbb{Z}_+$, $W_0$ is homeomorphic to $\overline{B}_n(0,1)$,
and (\ref{eq:nonAuto}) is $T$-periodic, then the orbit $v$ is
$T$-periodic.

\end{theorem}
\textbf{Proof:}
The existence of the full orbit in $W$ follows directly from Theorem~\ref{thm:discrete-trap-exists-orbit}. The existence of $T$-orbit
in $T$-periodic situation follows from Theorem~\ref{thm:discrete-traping-per}.

Let us denote by $\left[\frac{t}{T}\right]$ and $\left\{\frac{t}{T}\right\}$, the integer and fractional part of $\frac{t}{T}$. From Lemma~\ref{lem:lip1}
applied to $t_0=kT+ \left[\frac{t}{T}\right]T$ and $t_1=t_0 + \left\{\frac{t}{T}\right\}T$  and the estimate of the Lipschitz constants of $g_i$ we  obtain 
the following 
\begin{eqnarray*}
 \|\varphi(kT,t,z) - v(kT+t)\|=  \|\varphi(kT,t,z) - \varphi(kT,t,v(kT))\| = \\
  \left\| \varphi\left(kT,\left[\frac{t}{T}\right]T+ \left\{\frac{t}{T}\right\}T,z \right) -  
  \varphi\left(kT,\left[\frac{t}{T}\right]T+ \left\{\frac{t}{T}\right\}T,v(kT) \right) \right\| \leq \\
  \exp\left( \left\{\frac{t}{T}\right\}T B\right)   \left\| \varphi\left(kT,\left[\frac{t}{T}\right]T,z \right) -
  \varphi\left(kT,\left[\frac{t}{T}\right]T,v(kT) \right) \right\| \leq \\
   \exp\left( \left\{\frac{t}{T}\right\}T B\right) L^{\left[\frac{t}{T}\right]} \|z- v(kT)\|
\end{eqnarray*}

To obtain (\ref{eq:lip-iter}) from the above computations observe that
\begin{eqnarray*}
  L^{\left[\frac{t}{T}\right]}=\exp\left(\frac{\ln L}{T} \left[\frac{t}{T}\right] T   \right)= \\
  \exp \left(lt\right) \exp \left(-l \left(t - \left[\frac{t}{T}\right] T\right) \right).
\end{eqnarray*} 
\qed


\section{Self-consistent bounds and attracting orbits}
\label{sec:scb}

\subsection{The method of self-consistent bounds}
\label{subsec:method}

In this section we present an  adaption of the method of the
self-consistent bounds  \cite{ZM,Z2,Z3} to  non-autonomous
dissipative PDEs.

Let $J \subset \mathbb{R}$ be an interval (possibly unbounded). We
begin with an abstract nonlinear evolution equation in a real
Hilbert space $H$ (for example $L^2$) of the form
\begin{equation}
\label{eq:pde} \frac{du}{dt} = F(t,u),
\end{equation}
where the set of $x$ such that  $F(t,x)$ is defined for every $t
\in J$, denoted by $\widetilde{H}$, is  dense in $H$. Therefore
the domain
 of $F$ contains $J \times \widetilde{H}$.  By a solution of
(\ref{eq:pde}) we understand a function $u:J' \to \widetilde{H}$,
where $J' \subset J$ is an interval   such that $u$ is
differentiable and (\ref{eq:pde}) is satisfied for all $t \in J'$.

The scalar product in $H$  will be denoted by $(u | v)$.
Throughout the paper we assume that there is a set $I \subset
\mathbb{Z}^d$ and a sequence of subspaces $H_k \subset H$ for $k
\in I$, such that $\dim H_k \leq d_1 < \infty$ and $H_k$ and
$H_{k'}$ are mutually orthogonal for $k \neq k'$. Let $A_k: H \to
H_k$ be the orthogonal projection onto $H_k$. We assume that for
each $u \in H$ holds
\begin{equation}
   u=\sum_{k \in I} u_k=\sum_{k \in I} A_k u. \label{eq:H-decmp}
\end{equation}
The above equality  for a given $u \in H$ and $k \in I$ defines
$u_k$. Analogously if $B$ is a function with the range contained
in  $H$, then $B_k(u)=A_k B(u)$. Equation (\ref{eq:H-decmp})
implies that $H=\overline{\bigoplus_{k \in I} H_k}$.

Let  us fix an arbitrary norm on $\mathbb{Z}^d$, this norm will
denoted by $|k|$.

 For $n > 0$ we set
\begin{eqnarray*}
  X_n = \bigoplus_{|k| \leq n, k \in I } H_k, \\
  Y_n = X_n^\bot,
\end{eqnarray*}
by $P_n:H \to X_n$ and $Q_n:H\to Y_n$ we will denote the
orthogonal projections onto $X_n$ and onto $Y_n$, respectively.

\begin{definition}
\label{defn:Faddmissible} Let $J  \subset \mathbb{R}$ be an
interval. We say that $F: J \times H \supset \dom{F} \to H$ is
admissible, if the following conditions are satisfied for any $i
\in \mathbb{R}$, such that $\dim X_i >0$
\begin{itemize}
\item $ J \times X_i \subset \dom{F}$,
\item $P_i F : J \times X_i \to X_i$ is a $C^1$ function.
\end{itemize}
\end{definition}

\begin{definition}
Assume $F:J \times \widetilde{H} \to H$ is admissible. For a given
number $n>0$ the ordinary differential equation
\begin{equation}
  x'=P_n F(x), \qquad x \in X_n  \label{eq:galproj}
\end{equation}
will be called \emph{the $n$-th Galerkin projection} of
(\ref{eq:pde}).

By $\varphi^n(t_0,t,x_0)$ the solution of (\ref{eq:galproj}) with
the initial condition $x(t_0)=x_0$ at time $t_0+t$.
\end{definition}

\begin{definition}
\label{defn:selfconsistent} Assume $F:J \times \widetilde{H} \to
H$ is an admissible function. Let $m,M\in \mathbb{R}$ with $m\leq
M$. Consider an object consisting of:  $Z \subset J \times H,$  a
compact set $W\subset X_m$, such that $Z_t \subset W$ for $t \in
J$ and a sequence of compact sets $B_k \subset H_k$ for $|k|
> m$, $k \in I$. We define the conditions
\textbf{C1}, \textbf{C2}, \textbf{C3}, \textbf{C4a} as follows:
\begin{description}
\item[C1] For $|k|>M$, $k \in I$ holds $0 \in B_k$ .
\item[C2] Let $\hat{a}_k := \max_{a \in B_k} \|a\|$ for $|k| > m$, $k \in I$ and then
$\sum_{|k| > m, k \in I} \hat{a}_k^2 < \infty$. In particular
\begin{equation}
   W \oplus \Pi_{|k| > m} B_k \subset H
\end{equation}
and for every $u \in  W \oplus \Pi_{k \in I,|k| > m} B_k $ holds,
$\| Q_n u \| \leq \sum_{|k|>n, k \in I}\hat{a}_k^2$.
\item[C3] The function $(t,u)\mapsto F(t,u)$ is continuous on
$J \times W\oplus\prod_{k \in I, |k|>m}B_k \subset H$.

 Moreover, if we define for $k
\in I$, $f_k=\sup_{(t,u) \in J \times W \oplus \prod_{k\in I, |k|
>m} B_k} |F_k(t,u)| $, then $\sum f_k^2 < \infty$.

\item[C4] For $|k|>m$, $k \in I$  $B_k$ is given by (\ref{eq:C4aball}) or (\ref{eq:C4aint})
\begin{eqnarray}
  B_k&=&\overline{B(c_k,r_k)},  \quad r_k > 0  \label{eq:C4aball} \\
  B_k&=&\Pi_{s=1}^d[a_{s}^-,a_s^+], \qquad a_s^- < a_s^+, \:
  s=1,\dots,\dim(H_k) \label{eq:C4aint}
\end{eqnarray}
 Let $u \in W \oplus \Pi_{|k| > m}B_k$.
 Then for $|k| > m$ and $t \in J$ holds:
\begin{itemize}
\item if $B_k$ is given by (\ref{eq:C4aball}) then
   \begin{eqnarray}
     u_k \in  \partial_{H_k} B_k  & \Rightarrow  & (u_k -c_k | F_k(t,u)) <
     0. \label{eq:C4ain}
 \end{eqnarray}
 \item if $B_k$ is given by (\ref{eq:C4aint}), then for $t \in J$
 and  $s=1,\dots,\dim(H_k)$ holds
  \begin{eqnarray}
     u_{k,s}=a_{k,s}^-  & \Rightarrow  & F_{k,s}(t,u) >   0,  \label{eq:C4ain-} \\
     u_{k,s}=a_{k,s}^+  & \Rightarrow  & F_{k,s}(t,u) <   0.  \label{eq:C4ain+}
 \end{eqnarray}
 \end{itemize}
\end{description}
\end{definition}

In the sequel we will refer to equations (\ref{eq:C4ain}) and
(\ref{eq:C4ain-}--\ref{eq:C4ain+}) as  \emph{the isolation
equations} and to conditions \textbf{C1}, \textbf{C2}, \textbf{C3}
as \emph{the convergence conditions}.

Formally the above definitions require $Z \subset J \times X_m$,
but we will often apply them to $Z' \subset X_m$, so that we
assume that $Z=J \times Z'$ and the conditions C1,C2,C3,C4 refer
formally to the set $Z$.  In what follows quite often there will
be no need to distinguish these situations, and in such case we
will not  bother to state this explicitly, whether  $Z \subset
X_m$ or $Z \subset J \times X_m$.

 Given   $Z \subset J \times X_m$ (or $W \subset X_m$) and $\{B_k\}_{k \in I,
 |k|>m}$ satisfying conditions C1,C2,C3
 by $T$ (the tail) we will denote
\begin{equation*}
T:=\prod_{|k| > m} B_k \subset Y_m.
\end{equation*}

Here are some useful lemmas illustrating the implications of
conditions C1, C2, C3.

\begin{lemma}
\label{lem:compact} Let $W \subset X_m$ and $T \subset Y_m$. If
$W \oplus T$ satisfies condition \textbf{C2}, then $W\oplus T$ is
a compact subset of $H$.
\end{lemma}

\begin{lemma} \label{lem:compact2} Let $W \subset X_m$ and $T \subset
Y_m$.
Assume conditions C1,C2 and C3 on $W \oplus T$
for $F$ on $J$, then
\begin{displaymath}
  \lim_{n\to \infty} P_n(F(t,u)) =F(t,u), \quad \mbox{uniformly for
   $(t,u) \in J \times W \oplus T$}
\end{displaymath}
\end{lemma}

It turns out that for dissipative PDEs with periodic boundary
conditions it is rather easy to find $W \oplus T$ satisfying
C1,C2,C3,C4. We will have $e_k=\exp(ikx)$,
$\hat{a}_k=\frac{C}{|k|^s}$ and $\hat{f}_k=\frac{C}{|k|^{s-r}}$
with $s$ and $s-r$ as large as we want, to make the series $\sum_k
a_k \exp(ikx)$ converge uniformly together with some of its
derivatives. In particular, for $s$ sufficiently large $W\subset\subspaceH$ from
Theorem~\ref{thm:analyticTrappingRegion} forms self-consistent bounds, i.e.
it satisfies conditions C1,C2,C3,C4.

Observe that the topology on such set $W \oplus T$
for $s$ large enough is just the topology of the coordinate-wise
convergence. To be more precise we state the following lemma.

\begin{lemma}
\label{lem:scb-compact} Let $s>0$. Assume that $0 \leq \hat{a}_k
\leq \frac{C}{|k|^s}$ for $k \neq 0$ and $0 \leq \hat{a}_0 \leq
C$, and
\begin{equation}
   \sum_{k \in I} \hat{a}_k^2 < \infty.
\end{equation}
Let $Z=\{\{u_k\}_{k \in I} \ | \ |u_k| \leq \hat{a}_k \}$

 Then
 \begin{itemize}
 \item  $Z \subset H$, $Z$ is compact,
 \item  Let $\{z^j\}_{j \in \mathbb{N}} \subset Z$. Then $z^j \to
 z$ in $H$ iff for all $k \in I$  $z^j_k \to z_k.$
  \end{itemize}
\end{lemma}

\begin{definition}
\label{def:basic-diff-incl} Assume that $W \subset X_m$ and $T
\subset Y_m$. Let $W \oplus T$ satisfy conditions C1,C2,C3 for $F$
on $J$.

Let $c \in Y_m$, be such that  $c \in T$, $c=Q_m P_M c$ ( most of
the time we will take a center point of $T$).

Let $[\delta]=\{ P_m F(t,u + T) - P_mF(t,u + c)\ | \ u \in W  \}
\subset X_m $.

We define \emph{the basic differential inclusion} for
(\ref{eq:pde}) on $J \times W \oplus T$ as
\begin{equation}
  x'(t) \in P_m F(t,x(t)+c) + [\delta],  \quad x(t) \in  X_m. \label{eq:basic-diff-incl}
\end{equation}
and the translated $n$-th Galerkin projection of (\ref{eq:pde}) by
\begin{equation}
  x'(t)=P_n F(t,x(t)+Q_nc), \quad x(t) \in X_n.
  \label{eq:galproj-c-shifted}
\end{equation}
Let $\varphi_c^n$ be the local semiprocess induced by
(\ref{eq:galproj-c-shifted}).
\end{definition}

Observe that for $n >M$ holds $Q_n c=0$, hence
$\varphi_c^n=\varphi^n$.

The following two lemmas  clearly demonstrate the role of the
isolation condition C4. They show that it is enough to consider
the basic differential inclusion  (\ref{eq:basic-diff-incl}) to
build a trapping isolating segment (Lemma~\ref{lem:iso-seg}) or a
rigorous integrator (Lemma~\ref{lem:estm-proj-gal}). We omit
obvious proofs.

In our integration algorithm of dissipative PDE we compute
bounds for all Galerkin projections with $n>M$, hence we have in
fact $c=0$. Only when considering the $n$-th Galerkin projection
with $m<n<M$ we need to include $c$.

\begin{lemma}
\label{lem:iso-seg} Assume that $Z  \subset J \times X_m$, and
$Z\oplus T$ satisfies conditions C1,C2,C3 and C4.  Assume that $Z
$ is a trapping isolating segment for  differential inclusion
(\ref{eq:basic-diff-incl}).

Then for any $k>m$ the set $Z \oplus P_k T$ is a trapping
isolating segment for $\varphi_c^k$.
\end{lemma}

\begin{lemma}
\label{lem:estm-proj-gal} Assume that $W  \subset  X_m$ and
$W\oplus T$ satisfies conditions C1,C2,C3 and C4 on $J=[t_0,t_1]$.

Let $x_0 \in W$, be such that any $C^1$ solution of
(\ref{eq:basic-diff-incl}) with the initial condition $x(t_0)=x_0$
exists for $t \in [t_0,t_1]$ and is contained in $W$.

For $t \in [t_0,t_1]$ let
\begin{equation}
  x_I(t)=\{ y(t)\ | \ \mbox{$y$ is a $C^1$ solution of (\ref{eq:basic-diff-incl}), $y(t_0)=x_0$}\}
\end{equation}

Then for $k>m$ we have
\begin{equation}
  \varphi_c^n(t_0,t,x_0 + T) \in x_I(t_0+t) \oplus T, \quad t \in
  [0,t_1-t_0].
\end{equation}
\end{lemma}

Lemma~\ref{lem:estm-proj-gal} is the base on which our rigorous
integrator for dissipative PDEs with periodic boundary condition
is founded. For details how to estimate all solutions of
(\ref{eq:basic-diff-incl}) and how to estimate better the tail the
reader is referred to \cite{Z2,Z3,KZ}.

Sometimes, it will be  convenient to use a different norm on the
subspace containing $W \oplus T$ and $F(J \times (W \oplus T))$.
We just need to make sure that it induces the same topology on $W
\oplus T \cup F(J \times W \oplus T)$. This motivates the
following definition.
\begin{definition}
\label{def:compatiblenorms} Let $W \oplus T$ satisfy conditions
C1,C2,C3 for $F$ on $J$. We say that $\|x\|_1$ is a
\emph{compatible norm}  for $J \times (W \oplus T)$ if the
following conditions are satisfied
\begin{description}
\item[N1] $\|z\|_1$ is defined for $z \in W \oplus T \cup F(J \times (W \oplus T))$
\item[N2] $z^i \to z$ in $H$ for $z^i,z \in W \oplus T \cup F(J \times (W \oplus T)) $ iff $\|z_i - z\|_1 \to 0$
\item[N3] there exists $K$, such that for all $m\in\mathbb{R}$ holds $\|P_m z\|_1 \leq K
\|z\|_1$
\item[N4] $\|(I-P_n)\left( W \oplus T \cup F(J \times (W \oplus T))\right)\|_1 \to 0$ for $n \to 0$
\end{description}
\end{definition}

In the examples considered in our work we will have on $W \oplus
T$ the following estimate $|a_k| \leq \frac{C}{|k|^s}$ and on $F(J
\times (W \oplus S))$ the bounds  $|F_k| \leq
\frac{D}{|k|^{s-r}}$, where $s$ can be made as large as we want.
For example, the following are the compatible norms
$\|x\|=\sup_{k} |x_k|$ or $\|x\|=\sum_{k} |k|^p |x_k|$, for some
$p < s-r -d -1$, where $d$ is the dimension of the wave vectors
$k$ space.

\begin{lemma}
\label{lem:gal-proj-err-conv}  Let $W \subset X_m$. Assume that $W
\oplus T$ satisfies C1,C2,C3 for $F$ on $J$, $J$ is compact and
$\|\cdot \|_1$ is a compatible norm. Then
\begin{equation}
  \delta_n=  \sup_{(t,x) \in J\times (W \oplus T)} \|P_n(F(t,x)) - P_n(F(t,P_n
  x))\|_1 \to 0, \quad \mbox{for $n \to \infty$}
\end{equation}
\end{lemma}
\noindent \textbf{Proof:} From Lemma~\ref{lem:scb-compact} it
follows that on $W \oplus T \cup F(J \times W \oplus T)$ the
topology induced from $H$ coincides with the topology induced by
the norm $\|\cdot\|_1$. In the sequel we will use the distance
induced by this norm.

Observe that from condition \textbf{C3} it follows that $F:J
\times (W \oplus T) \to W \oplus T$ is continuous. Since $J
\times (W \oplus T)$ is compact, therefore $F$ on $J \times (W
\oplus T)$ is uniformly continuous, which expressed in terms of
the norm $\|\cdot\|_1$ means that for any $\epsilon>0$ there exits
$\delta(\epsilon)>0$, such that for any $t \in J$
\begin{equation}
\mbox{if }\|x-y\|_1 < \delta(\epsilon), \quad
 \mbox{then} \quad \|F(t,x)-F(t,y)\|_1 < \epsilon. \label{eq:equi-cont}
\end{equation}

Let us fix $\epsilon>0$. From condition \textbf{N4} it follows
that there exists $n_0$, such that for $n_0 >0$ the following
conditions are satisfied
\begin{eqnarray}
  \sup_{x \in W \oplus T} \| P_n(x) - x\|_1 < \delta(\epsilon),
  \quad n \geq n_0  \label{eq:Pn-conv} \\
   \sup_{ z \in F(J \times (W \oplus T))} \| P_n(z) - z\|_1 < \epsilon,
  \quad n \geq n_0. \label{eq:Pn-convF-2}
\end{eqnarray}
From (\ref{eq:equi-cont}) and (\ref{eq:Pn-conv}) it follows that
for any $t \in J$
\begin{equation}
   \sup_{x \in W \oplus T} \| F(t,P_n(x)) - F(t,x)\|_1 < \epsilon,
  \quad n \geq n_0 .
\end{equation}
This and (\ref{eq:Pn-convF-2}) imply that
\begin{equation}
 \sup_{(t,x) \in J\times (W \oplus T)} \|P_n(F(t,x)) - P_n(F(t,P_n
  x))\|_1 < 2 \epsilon, \quad n \geq n_0.
\end{equation}
 \qed

\subsection{Attracting orbits through trapping isolating segments}
\label{subsec:attr-orbits-iso-seg}

The goal of this section is to state the theorems that in
the context of self-consistent bounds and trapping isolating
segments will guarantee the existence of attracting orbit. The
orbit will be periodic if the forcing is periodic.

\begin{definition}
Consider (\ref{eq:pde}). Let $W \subset X_m$ and $W \oplus T
\subset H$ satisfy conditions C1,C2,C3 for $F$  on $J$. We say
condition {\em D} is satisfied on $J \times W \oplus T$ for the
compatible norm $\|\cdot\|_1$ if the following holds
\begin{description}
\item[D] There exists  $l \in \mathbb{R}$ such that for each Galerkin projection
\begin{equation}
  \sup_{(t,z) \in J \times W \oplus T}   \mu_1\left(\frac{\partial P_nF}{\partial z}(t,z):X_n \to X_n\right) \leq l
\end{equation}
\end{description}
where $\mu_1(A)$ for $A \in \mathbb{R}^{l\times l}$ is the
logarithmic norm of the matrix $A$ induced by the norm $\|\cdot
\|_1$.
\end{definition}

Condition \textbf{D} will be used to estimate the
Lipschitz constant of the semiflow induced by Galerkin projection
of our dissipative PDE and its Galerkin projection as discussed in
Section~\ref{subsec:LipConstFlow}. For this it is important that
set $W$ is convex.

\begin{theorem}
\label{thm:convergence} Let $W \subset X_m$ is convex, $T \subset
Y_m$ and $J=[t_0,t_1]$. Assume that on $J \times W \oplus T$
conditions C1, C2, C3 and condition $D$ are satisfied  for $F$ for
compatible norm $|\cdot|$.

Assume that for $n \in \mathbb{R}$ function $x^n:[t_0,t_1] \to  W
\oplus P_n T$ is a solution to $n$-th Galerkin projection of
(\ref{eq:pde}),
 such that $\lim_{n \to \infty} x^n(t_0) =x_0$.

 Then $x^n$
converge uniformly to $x:[t_0,t_1] \to  W \oplus T $, which is a
solution of (\ref{eq:pde}) and $x(t_0)=x_0$.
\end{theorem}
\noindent \textbf{Proof:} Let
\begin{displaymath}
\delta_n= \max_{(t,x) \in [t_0,t_1]\times W} |P_n(F(t,x)) -
P_n(F(t,P_n x))|.
\end{displaymath}
From Lemma~\ref{lem:gal-proj-err-conv} it follows that
\begin{equation}
 \lim_{n \to \infty}\delta_n= 0.
\end{equation}

Let us take $m \geq n$.  From Theorem \ref{thm:lognorm},  applied
to the $n$-th Galerkin projection of (\ref{eq:pde}) with $P_n x^m$
as 'an approximate solution' $y$,  it follows immediately that for
$t>t_0$ holds
\begin{equation}
  |x^n(t) - P_n(x^m(t))| \leq e^{l(t-t_0)}|x^n(t_0) - P_nx^m(t_0)| + \delta_n \frac{e^{l(t-t_0)} - 1}{l}
  \label{eq:estmdiffg}
\end{equation}

Observe that for $t \in [t_0,t_1]$ holds
\begin{eqnarray*}
 |x^n(t) - x^m(t)| \leq  |x^n(t) - P_n(x^m(t))| + |(I-P_n) x^m(t)| \leq \\
   e^{l(t-t_0)} |x^n(t_0) - P_n x^m(t_0)| +   \delta_n \frac{e^{l(t-t_0)} - 1}{l} + |(I-P_n) x^m(t)|
   \leq\\
      e^{l(t-t_0)} K |x^n(t_0) - x^m(t_0)| +   \delta_n \frac{e^{|l|(t_1-t_0)} - 1}{|l|} + |(I-P_n)(W\oplus T)|.
\end{eqnarray*}
This shows that $\{x^n\}$ is a Cauchy sequence in the norm $|
\cdot |$, hence it converges uniformly to $x:[t_0,t_1]\to W \oplus
T $. From Lemma 8 in \cite{ZNS} adopted to the non-autonomous
setting it follows that $\frac{d x(t)}{dt}=F(t,x(t))$. \qed

\begin{theorem}
\label{thm:LipschitzConstant}Let $J=[t_0,t_1]$ and   $Z \subset J
\times  X_m$ and $T \subset Y_m$, such that $Z_t$ is convex for $t
\in [t_0,t_1]$. Assume  that conditions C1, C2, C3 and condition
$D$ are satisfied on $Z \oplus T$ for a compatible norm $|\cdot|$.
Assume that $Z$ is a trapping isolating segment for
(\ref{eq:basic-diff-incl}).

Assume that functions $x:[t_0,t_1] \to  Z$ and $y:[t_0,t_1] \to Z$
are solutions of (\ref{eq:pde}).

Then for $t \in [t_0,t_1]$ holds
\begin{equation}
   |x(t) - y(t)| \leq e^{l(t-t_0)}|x(t_0) - y(t_0)|.
\end{equation}
\end{theorem}
\textbf{Proof:} From our assumption and Lemma~\ref{lem:iso-seg} it
follows that $Z \oplus P_n T$ for $n>M$ is a trapping isolating
segment for the $n$-th Galerkin projection of (\ref{eq:pde}).

 For $n >M$ let $x^n$ and $y^n$ be  solutions for
the $n$-th Galerkin projection of (\ref{eq:pde}) with the initial
conditions $x^n(t_0)=P_n x(t_0)$ and $y^n(t_0)=P_n y(t_0)$,
respectively. From Theorem~\ref{thm:lognorm} applied to the $n$-th
Galerkin projection with different initial conditions we obtain
\begin{equation}
  |x_n(t) - y_n(t)| \leq e^{l(t-t_0)}|P_n x(t_0) - P_n y(t_0)|.
  \label{eq:diffnlip}
\end{equation}
From Theorem~\ref{thm:convergence} it follows that  $x_n \to x$
and $y_n \to y$ uniformly on $[t_0,t_1]$. Then passing to the
limit in (\ref{eq:diffnlip}) gives
\begin{equation}
   |x(t) - y(t)| \leq e^{l(t-t_0)}|x(t_0) - y(t_0)|.
\end{equation}
\qed

\begin{theorem}
\label{thm:scb-exists-full-orbit} Let $J=\mathbb{R}$,   $Z \subset
J \times  X_m$ and $T \subset Y_m$, such that $Z_t$ is convex for
$t \in \mathbb{R}$. Assume  that conditions C1, C2, C3 and
condition $D$ are satisfied on $Z \oplus T$ for a compatible norm
$|\cdot|$. Assume that $Z$ is a trapping isolating segment for
(\ref{eq:basic-diff-incl}).

Then there exists $x:(-\infty,\infty) \to H$,   which is a
solution of (\ref{eq:pde}) and is contained in $Z \oplus T$, such
that for any solution $v:[t_0,\infty) \to H$ of (\ref{eq:pde})
with initial condition $v(t_0) \in Z_{t_0}$ holds
\begin{equation}
  |v(t_0+t)-x(t_0+t)| \leq e^{lt}|v(t_0)- x(t_0)|, \quad \mbox{for
  $t>0$} \label{eq:pde-attr-orb}
\end{equation}
where  $l$ is a constant bounding from above the logarithmic norm
in condition D.

In particular, if $l<0$, then the orbit $x$ attracts all solutions
in $Z \oplus T$.

If $Z$ is a $\Delta$-periodic and (\ref{eq:pde}) is
$\Delta$-periodic, then there exists $\Delta$-periodic orbit
contained $Z$.

\end{theorem}
\noindent \textbf{Proof:} From our assumption and
Lemma~\ref{lem:iso-seg} it follows that $Z \oplus P_nT$ for $n>M$
is a trapping isolating segment for the $n$-th Galerkin projection
of (\ref{eq:pde}). Therefore from
Theorem~\ref{thm:trap-exists-orbit} it follows that for any $n>M$
there exists $x^n:(-\infty,\infty) \to H$, such that $x(t) \in Z_t
\oplus T $ is a solution for the $n$-th Galerkin projection of
(\ref{eq:pde}).

From the Ascoli-Arzela lemma (compare the proof of
Theorem~\ref{thm:trap-exists-orbit}) it follows that the sequence
$\{x^n\}$ contains locally uniformly converging subsequence
 to $x:(-\infty,\infty) \to Z$. From
Lemma 8 in \cite{ZNS} adopted to the non-autonomous setting  it
follows that $\frac{d x(t)}{dt}=F(t,x(t))$.

Estimate \eqref{eq:pde-attr-orb} follows immediately from
Thm.~\ref{thm:LipschitzConstant}.

Observe $Z_0$ is homeomorphic to a closed finite-dimensional ball,
therefore the same is true for $Z_0 \oplus P_nT$. From this
observation  for $\Delta$-periodic trapping isolating segment and
$\Delta$-periodic equation from Theorem~\ref{thm:traping-per} we
obtain $\Delta$-periodic orbits $x^n:(-\infty,\infty) \to P_n(Z)$
for the $n$-th Galerkin projection of (\ref{eq:pde}). Now we apply
the Ascoli-Arzela lemma like in the first part of the proof. \qed

\subsection{Attracting orbits through  discrete time shifts}
 Assume that $N_0 \subset X_m$ is
compact, $N_0 \oplus T_0$, such that $|T_{0,k}| \leq
\frac{C_0}{|k|^{s_0}}$. $N\oplus T$ is our initial condition at
the time $t_0$.  One time step,  from $t=t_0$ to $t=t_0+h$,  of the rigorous
integrator described in \cite{Z2,Z3,C} does  the following
\begin{itemize}
\item[1.] Finds $W \subset X_m$ and $W \oplus T$, which satisfies conditions
C1,C2,C3,C4 and D for $F$ on interval $J=[t_0,t_0+h]$. Moreover, $N_0 \oplus
T_0 \subset W\oplus T$ and any solution of
(\ref{eq:basic-diff-incl}) with the initial condition $x(t_0) \in
N_0$ is defined for $t \in [t_0,t_0+h]$ and stays in $W$ for $t
\in [t_0,t_0+h]$.
\item[2.] From rigorous bounds for (\ref{eq:basic-diff-incl}) on
$J \times W \oplus T$ plus some linear uniform estimates for the
tail evolution, we obtain $N_1 \subset X_m$ and $T_1$, such that
$|T_{1,k}| \leq \frac{C_1}{|k|^{s_1}}$ and for any $n>M$ holds
\begin{equation}
  \varphi^n(t_0,h,N_0 \oplus P_n T_0) \subset N_1 \oplus P_n T_1
\end{equation}
\end{itemize}
It may happen that for a given $h>0$ the first stage might fail,
this part involves search for a priori bounds, which might not
exists if there is a blow-up for some solutions. This might happen
even for ODEs.

Therefore, our algorithm for rigorous integration of dissipative
PDEs, if completed with the success, give us uniform bounds for
solutions of all Galerkin projections.  Solutions for PDE satisfy
the same bounds as it follows from Theorem~\ref{thm:convergence}.
The same applies to the bounds for Lipschitz constants for the
semi-flow induced by the PDE and its Galerkin projections.

Now we will state the version of
Theorem \ref{thm:discrete-trap-attr-exists-orbit} for the context
of the method of self-consistent bounds

Let us fix $\Delta>0$. For any $n>M$ we define the discrete
semiprocess by setting
\begin{equation}
  g^n_i(x)=\varphi^n(i\Delta,\Delta,x), \label{eq:scb-T-time-shift}
\end{equation}
i.e. this a time shift by $t=\Delta$ from the section $t=i\Delta$
to $t=(i+1)\Delta$.

\begin{theorem}
\label{thm:scb-discrete-trap-attr-exists-orbit} Assume that there
exist compact and convex set $Z \subset X_m$  and $T\subset Y_m$,
such that conditions C1,C2,C3 and D for some compatible norm
$\|\cdot\|$ are satisfied  on $Z \oplus T$ for F on
$J=\mathbb{R}$.

 Assume $W_i \subset X_m$ and $T_i \subset Y_m$ for
 $i=0,1,\dots,k-1$ are such that for all $i=0,\dots,k-1$ and $n>M$ holds
 \begin{eqnarray}
   W_i \subset Z, \quad T_i \subset T  \\
   \varphi^n(i\Delta,[0,\Delta],W_i \oplus P_n T_i) \subset Z
   \oplus T,  \\
   g^n_i(W_i \oplus P_n T_i) \subset W_{(i+1)\mod k} \oplus P_n
   T_{(i+1)\mod k}.
 \end{eqnarray}

Let $L_i,B \in
\mathbb{R}$ be such that for $i=0,1,\dots,k-1$ holds
\begin{eqnarray}
  \sup_{x \in W_i} \|Dg_i(x)\| &\leq& L_i , \\
   \sup_{(t,z) \in \mathbb{R} \times
  Z \oplus T} \mu\left(\frac{\partial P_n F}{\partial
z}(t,P_n z): X_n \to X_n\right) &\leq& B.
\end{eqnarray}

Then there  exists a full orbit $v$ for $\varphi$,
$C=\max(1,\exp(B \Delta))$ and $l=\frac{1}{\Delta}\ln ( L_0 L_1
\dots L_{k-1})$, such that for any $i=0,\dots,k-1$  and $ z \in
W_i\oplus T_i$ and $t
>0$
\begin{equation}
 \|\varphi(i\Delta,t,z) - v(i\Delta+t)\| \leq C\exp (l t) \|v(i\Delta) -  z\|,
\end{equation}
where $\varphi$ denotes the semiprocess induced by (\ref{eq:pde}).

If $l<0$ the orbit $v$ attracts all orbits starting from  $W_i
\oplus T_i$.

If $W_0$ is homeomorphic to $\overline{B}_n(0,1)$
and (\ref{eq:pde}) is $\Delta$-periodic, then the orbit $v$ is
$\Delta$-periodic.

\end{theorem}


\section{Proof of Theorem~\ref{thm:main1}}
\label{sec:proof}

The proof follows the scheme of the proof of \cite[Theorem
1.1]{C}.  The important modification is the inclusion of
non-autonomous forcing, which requires the estimates for the
Lipschitz constant for the flow discussed in
Section~\ref{sec:Lip-const-flow} and
Section~\ref{subsec:attr-orbits-iso-seg}. We will argue here that
all the computer assisted checks need to obtain
Theorem~\ref{thm:main1} which is stronger in conclusions  than
\cite[Theorem 1.1]{C} are already contained in the proof of
\cite[Theorem 1.1]{C}.

\providecommand{\gp}{Galerkin projection of \eqref{eq:burgers_infinite1}}

\paragraph{Proof}

The three main steps in the proof  are as follows\\

\begin{enumerate}
  \item Construction of an absorbing set, $\mathcal{A}\subset\subspaceH$, see
  Definition~\ref{def:absorbingSet}.
  \item Construction a time independent trapping isolating segment
  $W \subset \subspaceH$ and establishing the existence of the locally
  attracting orbit within $W$. For this  we use
  Theorem~\ref{thm:scb-exists-full-orbit} and we check
  that $l<0$ (this is the bound for the logarithmic norms).
  \item Rigorous numerical integration of the  absorbing set $\mathcal{A}$ up to the time
    when interval bounds for the solutions of the partial
    differential equation are contained in the interior of trapping isolating segment $W$.
\end{enumerate}
In what follows we discuss the above three steps separately.

\paragraph{Step 1} The existence an absorbing set $\mathcal{A}$ is established in
Theorem~\ref{thm:absorbingSet} and an algorithm for its
construction is presented in \cite[Section~8]{C}. The only
constants depending on the forcing appearing in the construction
of $\mathcal{A}$ are the energy of the forcing
($E_0=E(\{f_k\})\nu^{-2}$) -- in assumptions of
Theorem~\ref{thm:absorbingSet}, and the absolute value of the
forcing modes $|f_k|$ appearing in \eqref{eq:depOnF1}, and
\eqref{eq:depOnF2}.

In the problem \eqref{eq:burgers_infinite1} we split the forcing
into the autonomous part $f(x)$, and the nonautonomous part
$\widetilde{f}(t,x)$, such that $f(t,x)=f(x)+\widetilde{f}(t,x)$.
According to our assumptions we have
\begin{equation}
\widetilde{f}_k(t)\in[-\varepsilon, \varepsilon], \quad
\varepsilon=0.03, \qquad \forall k\in\mathbb{Z}, \ t \in
\timeInterval,
\end{equation}
Thus the following constants required in the
construction can be easily bounded
\begin{itemize}
  \item total energy of the forcing $E(\{f_k + \widetilde{f}_k(t)\})\leq E(\{|f_k|+\varepsilon\})$ for all $t\in\timeInterval$,
  \item absolute value of the forcing contribution to $\frac{d a_k}{d t}$: $|f_k + \widetilde{f}_k(t)|\leq |f_k|+\varepsilon$, for all
    $t\in[0,\infty),\ k\in\mathbb{Z}$.
\end{itemize}
Having these bounds, the algorithm from Section~8 in \cite{C} is applied directly.

\paragraph{Step 2}
Construction of the trapping isolating segment, $W$. This involves
verifying that the vector field points inwards on the boundary of
the trapping isolating segment. The trapping isolating segment is
required to be of the form $W=\mathbb{R} \times W_0$.
Observe that the
right-hand side of \eqref{eq:burgers_infinite1} has to be
evaluated for all times $t\in[0,\infty)$. This is achieved by
using the interval arithmetic, and plugging-in the interval bound
$[-\varepsilon,\varepsilon]$ in place of $\widetilde{f}_k(t)$ for
all $k\in\mathbb{Z}$, thus the obtained set is time-independent.
The attraction toward the
fixed point is obtained by the computation of logarithmic norm $l$. If $l<0$, then we 
 just apply Theorem~\ref{thm:scb-exists-full-orbit}.

\paragraph{Step 3}
For the rigorous numerical integration we have been using
\emph{the Lohner-type algorithm for differential inclusions}
proposed in \cite{KZ,Z3}. The differential inclusion is needed to treat the nonautonomous part for which
we just have the bound $\widetilde{f}_k(t)\in[-\varepsilon, \varepsilon]$.
 In \cite{KZ} it is argued that this
algorithm works for time dependent perturbations for which there
is an a-priori knowledge that they can be contained in an interval
box.

In the algorithm for rigorous numerical
integration of dPDEs that we used \cite[Algorithm~1]{C}, the
contribution of the nonautonomous forcing is accordingly added to
the actual perturbations vector, see Step~5 of
\cite[Algorithm~1]{C}.

\qed

An interesting consequence of the fact that the computer assisted
part of proof of Theorem~\ref{thm:main1} is essentially the same
as for the proof of \cite[Theorem 1.1]{C} is that  all example
theorems from \cite{C}, presented in the table \cite[Table 1]{C}
are true for a much wider class of forcing functions than it was
claimed in \cite{C}, but we have to replace the fixed
point by the periodic orbit for the time-periodic forcing and
simply attracting orbit for the non-periodic forcing. Namely,
they are true for the nonautonomous forcing, consisting from
autonomous and non-autonomous parts. In \cite{C} the nonautonomous
part satisfied $|\tilde{f}_k| \leq \varepsilon$ for $0 \neq |k|
\leq m$.  The values of $\varepsilon$ are provided in the table
\cite[Table 1]{C} for each example theorem that was proved.


\providecommand{\mappingSymbol}{\Phi_{t_p}}
\providecommand{\mappingSymboln}[1]{\Phi_{#1,t_p}}
\providecommand{\imageTail}{T_{\mappingSymbol}}
\providecommand{\rhsSymbol}{G}
\providecommand{\trappingRegion}{W}
\providecommand{\tubeSymbol}{W}
\providecommand{\absorbingSet}{\mathcal{A}}
\providecommand{\firstSegment}{\trappingRegion_0}
\providecommand{\family}{\mappingSymboln{jt_p},\ j \in \mathbb{Z}}
\providecommand{\familyShort}{\mappingSymboln{jt_p}}

\section{Algorithm for the proof of  Theorem~\ref{thm:main2}}
\label{sec:proof-main2}

\begin{definition}
  \label{def:Pmap}

  Let $t_0 \in \mathbb{R}$, $t_p> 0$, $x\in\subspaceH$. According with the notation introduced in
  Section~\ref{subsec:semi-proc} by
  $\varphi(t_0,t_p, x)\in\subspaceH$ we denote the
  time shift by $t_p$ along the solution of \eqref{eq:burgers_infinite1}
   with i.c. $x(t_0)=x_0\in\subspaceH$, which is defined due to the existence and the uniqueness
   of solutions of \eqref{eq:burgers_infinite1} within the
  subspace $\subspaceH$.

  We define $\mappingSymboln{t_0}\colon\subspaceH\to\subspaceH$ as
  \begin{equation}
    \mappingSymboln{t_0}\colon x\mapsto\varphi(t_0,t_p, x).
  \end{equation}
\end{definition}

The proof of Theorem~\ref{thm:main2} will have the same three step
structure just as the proof of Theorem~\ref{thm:main1}. However, in the
case of Theorem~\ref{thm:main2} we have the time dependent forcing
term, which cannot treated as a small perturbation of the
autonomous part of (\ref{eq:burgers_infinite1}).
 The main difference is that now we
will consider the family of maps $\family$
the time shift by the period of the forcing (or the period of the main part of
the forcing).

  For the family of maps $\family$ we establish the
existence of the absorbing set $\mathcal{A}$ (step 1), the
existence of a trapping isolating segment
$\trappingRegion\subset\mathbb{R}\times\subspaceH$ in which the family of maps
$\family$ is a contraction (step 2) and we show that
$\mappingSymboln{(j+n-1)t_p}\circ\cdots\circ\mappingSymboln{(j+1)t_p}\circ\mappingSymboln{jt_p}(\mathcal{A})\subset W_0$ for a $n \geq 1$, 
and all $l\in\mathbb{Z}$
(step 3). We denote 
\begin{equation*}
  \mappingSymbol^n(\mathcal{A}):=\mappingSymboln{(j+n-1)t_p}\circ\cdots\circ\mappingSymboln{(j+1)t_p}\circ\mappingSymboln{jt_p}(\mathcal{A}).
\end{equation*}
This is a little abuse of notation, but we hope that it will not cause any misunderstanding. Observe that in this case in order to calculate the
Lipschitz constant of $\family$ we cover an approximate time
periodic solution of the problem \eqref{eq:burgers_infinite1} by a
finite number of interval enclosures, and then estimate the
logarithmic norms locally for each piece.
 In such setting Step 1 is the same in both proofs. Here Step 2 requires
 the computation of the uniform bounds for $\mappingSymboln{jt_p}$ for $ j \in \mathbb{Z} $ and its Lipschitz constant
 (as in \cite{C}) we use the logarithmic norms for that. The
 computation of $\family$ is done with our rigorous integrator
 for dPDEs. In Step 3 we again do the rigorous integration of dPDEs
 to compute $\family$.

In order to obtain rigorous bounds for the family of maps $\family$ a $C^0$
rigorous numerical integrator, capable of integrating
nonautonomous system of equations, has to be employed.   This can be
achieved by \cite[Algorithm~1]{C} with just one modification.
Instead of using \emph{the $C^0$ Lohner integrator} in
\cite[Algorithm~1, Step~4]{C} to solve the system of autonomous ODEs,
a $C^0$ Lohner nonautonomous integrator is used to solve the system of nonautonomous ODEs.
Technically, the automatic differentiation of the nonautonomous forcing
modes is performed in this step in order to calculate the higher order time derivatives,
and include the contribution of the nonautonomous forcing into the Taylor coefficients.

Below, we present an algorithm for proving
Theorem~\ref{thm:main2}.

\paragraph{\textit{Notation}} Let $0<m\leq M$.
Following the notation from Section~\ref{sec:scb} by $W\subset\mathbb{R}\times\subspaceH$
we denote a trapping isolating segment, and
$W_0=\left\{x\in\subspaceH\colon (0,x)\in W\right\}$.
By $[W]\subset\mathbb{R}\times\subspaceH$ we denote a representation
of $W$ in the algorithm (interval bounds enclosing the trapping
isolating segment $W$).
We assume that $W_0$ forms self-consistent
bounds, and can be divided into finite part $P_mW_0\subset
P_m\subspaceH$, and the infinite dimensional part (the tail)
$T_{W_0}:=(I-P_m)W_0\subset (I-P_m)\subspaceH$.
Similar to our previous works, for
technical reasons, the finite part of the tail (indexed by
$i\colon m<i\leq M$) is distinguished from the infinite
dimensional part. For the definition of self-consistent bounds
refer to Section~\ref{sec:scb}. Let $i\in\mathbb{Z}$.
By $(W_0)_i\subset\mathbb{R}^2$ we denote the $i$-th coordinate of $W_0$.
Although the coordinates of $W_0$ are pairs of real numbers representing
the real and imaginary parts of complex numbers, often we
use the notation $(W_0)_i^{+/-}$, meaning that the
corresponding operations are performed for the real and the
imaginary parts separately, and $(W_0)_i^{+/-}$ returns
supremum/infimum of the real and the imaginary part of $(W_0)_i$ respectively.

In the algorithm description, to simplify the notation, we will drop
the first part of subscript in
the symbol denoting the family of maps $\familyShort$, and use simply
$\mappingSymbol$, as we are performing the
rigorous numerical integration for all initial times $jt_p$ simultaneously.

We will use the notation $[\mappingSymbol(\cdot)]$ to denote
rigorous interval bounds for the image of $\mappingSymbol$
obtained by applying the rigorous $C^0$ Lohner integrator,
$P_m\mappingSymbol$ denotes the finite dimensional version of
$\mappingSymbol$, in the sense that $m$-th Galerkin projection of
\eqref{eq:burgers_infinite1} is integrated in order to calculate
the image. We denote the tail of
$\left[\mappingSymbol(W_0)\right]$ outputted from our rigorous
integrator by
$\imageTail:=(I-P_m)\left[\mappingSymbol(W_0)\right]$,
$C_T>0$, and $s_T>0$ denote the constants defining the polynomial
bound for the tail $T$, our algorithm is able to calculate efficiently
this values
(refer technical description in \cite[Appendix B and Appendix C]{C}).
inflate$\left((W_0)_i, c_3\right)$
inflates $(W_0)_i$ -- $i$-th coordinate (a complex number)
of $W_0$, i.e. makes it wider by the constant $c_3>0$.
By $\rhsSymbol\colon\mathbb{R}\times\subspaceH\to\subspaceH$
we will denote the right hand side of
\eqref{eq:burgers_infinite1}. Following the notation from Section~\ref{sec:Lip-const-flow}
by $\mu(A)$ we denote the
logarithmic norm of a square matrix $A$, $\mu_{b,\infty}$ is the
logarithmic norm inducted by the so-called block-infinity norm,
see \cite{ZAKS} for details.

\begin{figure}[htbp]
  \centering
  \includegraphics[width=\textwidth]{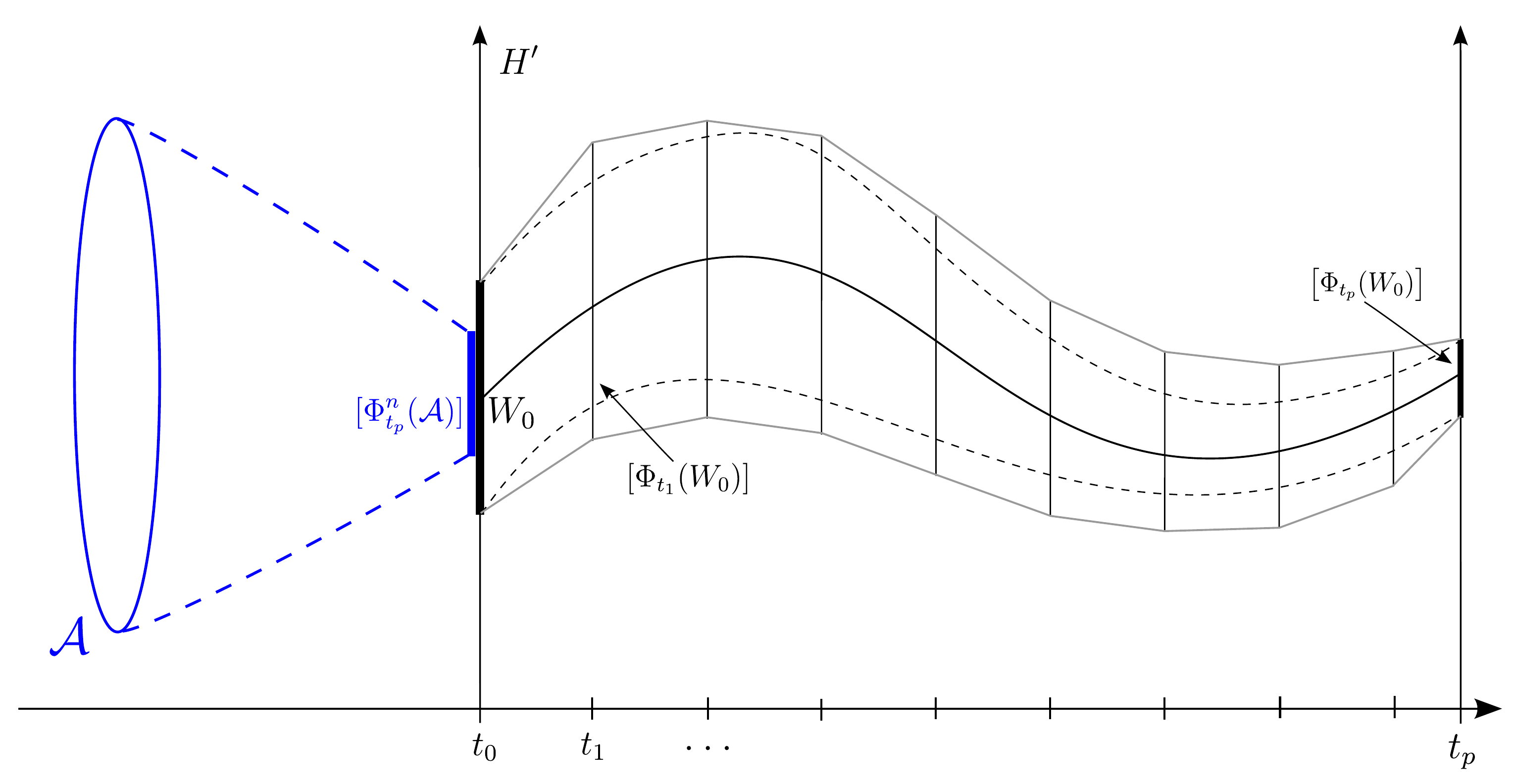}
  \caption{Diagram presenting steps of Algorithm for the proof of  Theorem~\ref{thm:main2}}
  \label{fig:first}
\end{figure}	

\subsection{Main algorithm}
\label{sec:algMainAlgorithm}

\paragraph{Input}
    \begin{itemize}
        \item $M\geq m>0$, integers, $m$ -- the Galerkin projection \eqref{eq:symmetricGalerkinProjection} dimension, and
          $M$ -- the dimension of the finite tail part of self-consistent bounds,
        \item $[\nu_1,\nu_2]>0$, an interval of the viscosity constant values -- can be degenerate (single valued),
        \item $s\geq 4$, the \decay\ of coefficients that is required from the constructed bounds and trapping regions,
        \item order and the time step of the Taylor method used by the \emph{$C^0$ Lohner nonautonomous integrator},
        \item $t_p>0$ period of the nonautonomous forcing,
        \item the forcing $F$ modes, the autonomous part is provided by $[f_k]:=f_k+[f_\varepsilon]$, where
        $[f_\varepsilon]$ is a uniform and constant perturbation $[f_\varepsilon]=[-\varepsilon, \varepsilon]\times[-\varepsilon,\varepsilon]$,
        and the nonautonomous part is provided by finite number of nonautonomous sufficiently regular $t_p$-periodic in time forcing modes
        $\{\widetilde{f}_k(t)\}_{0<|k|\leq m}$, given explicitly in a closed, representable on computer form, allowing automatic differentiation,      
      \item the constants $c_1>0$, $c_2>0$, and $c_3>0$, which should be adjusted according to the equation considered (in our program we have
        $c_1=10^{-5}$, $c_2=0.1$, $c_3=1.01$).
   \end{itemize}
\paragraph{Output}
    \begin{itemize}
        \item $[\trappingRegion]\subset[0,t_p]\times\subspaceH$ -- interval representation of $W\subset\mathbb{R}\times\subspaceH$
        -- the trapping isolating segment for the discrete semi-process, by representation
        of a trapping isolating segment we mean interval bounds enclosing $W$, and enclosing all trajectories traversing the
        trapping isolating segment, this set when glued together bound $W$, $[\trappingRegion]$ is used to calculate the
        Lipschitz constant bounds,
        \item $\mathcal{A}\subset\subspaceH$, an absorbing set forming self-consistent bounds for \eqref{eq:burgers_infinite1},
        \item a upper bound for the Lipschitz constant $L=e^l$ of $\mappingSymbol$ on the set $\firstSegment$,
    \end{itemize}
  \paragraph{begin}
    \begin{enumerate}
        \item
         Using the algorithm from Section~\ref{sec:algTrappingRegion} calculate a set $\firstSegment\subset\subspaceH$
        in which $\mappingSymbol$ is contraction. Verify that $L=e^l$ the calculated bound for the Lipschitz constant of
        $\mappingSymbol$ on $\firstSegment$ satisfies the inequality $L<1$, and therefore there is an (locally) attracting periodic
        solution within the trapping isolating segment $\trappingRegion$.

        \item Using the procedure from \cite[Section~8]{C} calculate the absorbing set
        $\absorbingSet\subset\subspaceH$ proved to exist in Thm.~\ref{thm:absorbingSet},
        taking care of the nonautonomous part of forcing. Put
        \begin{equation}
          \label{eq:absorbingBallRadius}
          E_0=\sup_{\substack{g\in\left\{[f_k]\right\}\\t\in\mathbb{R}}}{\frac{E(\{g_k+\widetilde{f}_k(t)\})}{\nu^2}},
        \end{equation}
        and wherever value of $\left|[f_k]\right|$ is required, put $\left|[f_k]+\sup_{t\in\mathbb{R}}{\widetilde{f}_k(t)}\right|$.

        \item Using the \emph{$C^0$ Lohner nonautonomous integrator} calculate\\$[\mappingSymbol^n(\absorbingSet)]$
        until $n>0$ is found such that $[\mappingSymbol^n(\absorbingSet)]\subset \firstSegment$.
    \end{enumerate}
  \paragraph{end}

\paragraph{}According to Theorem~\ref{thm:absorbingSet} procedure from \cite[Section~8]{C} can be modified -- value of $\alpha$ \eqref{eq:fixedInt}
can be omitted in the estimates, but then according to Lemma~\ref{lem:akbk2} there is a penalty -- the estimates for the norm are
received (compare \eqref{eq:akt-exp-small}) instead of estimates for the infimum and supremum.

For cases with $\alpha$ small, the original procedure from \cite[Section~8]{C} is expected to be more efficient. Whereas, for cases
with $\alpha$ large the estimates based on Lemma~\ref{lem:akbk2} are expected to be more efficient. The recommended strategy is to
calculate estimates using both of the above presented methods, and then take the intersection.

Steps of the algorithm described above are graphically presented on Figure~\ref{fig:first}.

\subsection{Algorithm constructing bounds for trapping isolating segment  and estimating  Lipschitz constant for time shifts}
\label{sec:algTrappingRegion}

\paragraph{Input}
    The same as in the Algorithm from Section~\ref{sec:algMainAlgorithm}.
\paragraph{Output}
    \begin{itemize}
        \item $[\trappingRegion]\subset[0,t_p]\times\subspaceH$ -- interval representation of a trapping isolating segment
        for the discrete semi-process,
        \item $L>0$ a bound for the Lipschitz constant of $\mappingSymbol$ on $\firstSegment$,
    \end{itemize}

    \paragraph{begin}
    \begin{enumerate}
      \item Find $\overline{x}$ an approximate location of the fixed point of $\mappingSymbol$ by applying
      the Newton method to the map
        $P_m\mappingSymbol(x)-x$, i.e.
        \begin{displaymath}
        x^{(k+1)}=x^{(k)}- (DP_m\mappingSymbol\left(x^{(k)}\right)-Id)^{-1}\left(P_m\mappingSymbol\left(x^{(k)}\right)-x^{(k)}\right),
        \end{displaymath}
         stop after
        several iterations.\\

         \item\label{step:mappedItself} Iteratively find a set
         $\firstSegment\subset\subspaceH$, such that
          $\left[\mappingSymbol(\firstSegment)\right]\subset~\firstSegment$
          using the following procedure.

          As the initial value of $\firstSegment$ take an interval hull of $\overline{x}$
          in $\subspaceH$. Initialize $\firstSegment$ by adding
          $[-c_1, c_1]\times[-c_1, c_1]$ to all coordinates of $\overline{x}$,
          initialize the tail part $T:=T_{\firstSegment}$ with values such that $T_i$ satisfy
          $|T_i|\leq C/|i|^s$, where $C:=|\left(\firstSegment\right)_m|\cdot m^s$ for $i>m$, and $s$ is
          provided as an input, $\left(\firstSegment\right)_m$ denotes the $m$-th coordinate of $\firstSegment$.

        \textbf{while} $\left[\mappingSymbol(\firstSegment)\right]\nsubseteq \firstSegment$\textbf{ do }
        \par
    \begingroup
    \leftskip2em
    \rightskip\leftskip
      \textbf{for each} $i\in\{1,\dots,M\}$\textbf{ such that }$\left[\mappingSymbol(\firstSegment)\right]_i\nsubseteq(\firstSegment)_i$\textbf{ do }\\
      \par
      \begingroup
      \leftskip4em
      \rightskip\leftskip
        \textbf{if }$\left[\mappingSymbol(\firstSegment)\right]_i^-\leq(\firstSegment)_i^-$\textbf{ then }\\
        \par
        \begingroup
        \leftskip6em
        \rightskip\leftskip
        $(\firstSegment)_i^-:=\left[\mappingSymbol(\firstSegment)\right]_i^-+c_2\cdot\left(\left[\mappingSymbol(\firstSegment)\right]_i^--(\firstSegment)_i^-\right)$\\
        inflate$\left((\firstSegment)_i, c_3\right)$\\
        \par
        \endgroup
        \textbf{end}\\
        \textbf{if }$\left[\mappingSymbol(\firstSegment)\right]_i^+\geq(\firstSegment)_i^+$\textbf{ then }\\
        \par
        \begingroup
        \leftskip6em
        \rightskip\leftskip
        $(\firstSegment)_i^+:=\left[\mappingSymbol(\firstSegment)\right]_i^++c_2\cdot\left(\left[\mappingSymbol(\firstSegment)\right]_i^+-(\firstSegment)_i^+\right)$\\
        inflate$\left((\firstSegment)_i, c_3\right)$\\
        \par
        \endgroup
        \textbf{end}\\
      \par
      \endgroup
      \textbf{end}\\
      \textbf{if }$C_T\cdot(M+1)^{-s_T}\leq C_{T_{\mappingSymbol}}\cdot(M+1)^{-s_{T_{\mappingSymbol}}}$\textbf{ then }\\
      \par
      \begingroup
      \leftskip6em
      \rightskip\leftskip
      $C_T:=C_{T_{\mappingSymbol}}\cdot(M+1)^{s_T-s_{T_{\mappingSymbol}}}$\\
      \par
      \endgroup
      \textbf{end}
    \par
    \endgroup
        \textbf{end}.\\

        \item Using the \emph{$C^0$ Lohner nonautonomous integrator} rigorously integrate $\firstSegment$ to obtain bounds along the
        orbit for the times $t_1<\dots<t_n<t_{n+1}=t_p$, which we denote by $\left[x_{i}\right]$ for $i=1,\dots,n+1$, and
        $\left[x_0\right]:=W_0$.

        As the output of our Lohner type algorithm for dPDEs the so-called \emph{rough enclosure} -- rigorous
        bounds for the solution during the whole time step are obtained, we denote the obtained rough-enclosures by
        $\left[\varphi\left(t_i,[0, t_{i+1}-t_i],[x_{i}]\right)\right]$ for $i=0,\dots,n$. For $i=0,\dots,n$ the set
        $\left[\varphi\left(t_i,[0, t_{i+1}-t_i],[x_{i}]\right)\right]$ forms self-consistent bounds
        for $G$ on $[t_i, t_{i+1}]$, i.e. satisfy conditions
        C1, C2, C3 from Definition~\ref{defn:selfconsistent}.

        In order to calculate the Lipschitz constant of $\mappingSymbol$ on $\firstSegment$ we construct the following interval bounds
        enclosing the trapping isolating segment for the discrete semi-process, and enclosing all trajectories traversing the
        trapping isolating segment.

        \begin{equation*}
          [\tubeSymbol]:=\bigcup_{i=0}^n[t_i, t_{i+1}]\times\left[\varphi\left(t_i,[0, t_{i+1}-t_i],[x_{i}]\right)\right].
        \end{equation*}
        Then the logarithmic norms are calculated locally on each part of $[\tubeSymbol]$.
        \item Using the bounds calculated in the previous step calculate local logarithmic norms in a suitable norm.
        First, calculate an orthogonal change of coordinates $Q_0$ such that
        \begin{equation}
          \label{eq:diagonalForm}
          [Q_0^{-1}]\cdot \frac{\partial P_mG}{\partial a}\left(\cdot, \left[\varphi\left(0,[0, t_1],\firstSegment\right)\right]\right)\cdot Q_0
        \end{equation}
        is in the close to the block-diagonal form. In the equation \eqref{eq:diagonalForm} $Q_0$ denotes a point matrix composed of approximate
        eigenvectors, possibly obtained from a non-rigorous external numerical package, $[Q_0^{-1}]$ is the rigorous (interval) inverse of $Q_0$,
        we put $\cdot$ as the first argument of $\frac{\partial P_mG}{\partial a}$, because this value is irrelevant (our assumption that the
        nonautonomous term does not depend on $a$ at all).\\

        \textbf{for each} $i\in\left\{0,\dots,n\right\}$ \textbf{ do }\\
    \begin{equation}
      \label{eq:localLogNorm}
      \quad\text{calculate }l_{i}:=\mu_{b,\infty}\left(\frac{\partial G}{\partial a}\left(\cdot, \left[\varphi\left(t_i, [0, t_{i+1}-t_i],\left[x_i\right]\right)\right]\right)\right)
    \end{equation}
    \textbf{end},\\
    where $\mu_{b,\infty}$ is the logarithmic norm inducted by the block infinity norm defined using the orthogonal change of coordinates $Q_0$
    (that norm is denoted in Lemma~\ref{lem:lip1} by $||\cdot||_0$). Obviously $||\cdot||_0$ is a compatible norm according to
    Definition~\ref{def:compatiblenorms}.

    Observe that at a step $j\colon0<j\leq n$ in the integration process the matrix $Q_0$ can be replaced with another matrix , such that
    the matrix
    \begin{equation}
      \label{eq:diagonalForm2}
      [Q_j^{-1}]\cdot \frac{\partial P_mG}{\partial a}\left(\cdot, \left[\varphi\left(t_j,[0, t_{j+1}-t_j],\firstSegment\right)\right]\right)\cdot Q_j
    \end{equation}
    is in the close to a block-diagonal form. Observe that in this case $||\cdot||_0=\dots=||\cdot||_{j-1}$, and the local logarithmic
    norms $l_{j-1}$ and $l_j$ are calculated using two distinct norms -- $||\cdot||_0$ and $||\cdot||_j$ respectively.
    \item Calculate the global Lipschitz constant using the local logarithmic norms calculated in the previous step. Depending on the
      number of distinct norms that were used to calculate the local logarithmic norms, two cases are distinguished.

      \emph{Case I} -- only the norm $||\cdot||_0$ was used to calculate all of the logarithmic norms $\{l_{i}\}_{i=0}^{n}$ in the previous
      step. According to Theorem~\ref{thm:scb-discrete-trap-attr-exists-orbit} the Lipschitz constant of $\mappingSymbol$ is bounded by 
      $L=Ce^{\frac{l}{\Delta}\Delta}=Ce^{l}$,
      where $\Delta=t_p$, $l=\sum_{i=0}^n{l_{i}\cdot (t_{i+1}-t_i)}$, and $C\geq 1$.\\

      \emph{Case II} -- at least two norms were used to calculate the logarithmic norms $\{l_{i}\}_{i=0}^{n}$ in the previous step of the
      algorithm. According to Theorem~\ref{thm:scb-discrete-trap-attr-exists-orbit} the Lipschitz constant of $\mappingSymbol$ is bounded by
      $L=Ce^{\frac{l}{\Delta}\Delta}=Ce^{l}$, where $\Delta=t_p$, $l=\sum_{i=0}^n{l_{i}\cdot (t_{i+1}-t_i)P_{i\mapsto i+1}}$, and $C\geq 1$.\\

      If the norms $||\cdot||_j$ and $||\cdot||_{j+1}$ are different put
      $P_{j\mapsto j+1}:=||Q_{j+1}^{-1}Q_j||_\infty$, otherwise, put $P_{j\mapsto j+1}:=1$.\\

      If $l<0$ then the existence of a locally attracting orbit within the set $\tubeSymbol$ is claimed.

    \end{enumerate}
\paragraph{end}

\begin{rem}
            All the bounds for the logarithmic norm of the (infinite dimensional) derivative of vector field $D\rhsSymbol$ calculated
            in the main algorithm presented above are carried out in suitable block coordinates. The finite part of $DP_mG$ is reduced
            by an orthogonal change of coordinates to an (almost) block-diagonal form, i.e. having $2\times 2$ blocks on the
            diagonal. The block decomposition of $H$ is given by $H=\oplus_{(i)}H_{(i)}$.
            For $(i)\leq m$
            each block $H_{(i)}$ is a two-dimensional eigenspace of $J$. In case of two dimensional blocks $(i)=(i_1,i_2)\in\mathbb{Z}^2$,
            the expression $(i)<m$ means that $i_j<m$ for $j=1,2$. We consider all blocks two dimensional, and for $(i)>m$
            (the infinite dimensional part) the diagonal
            blocks look like
            $\left[\begin{array}{cc}\lambda_i&\alpha_i\\-\alpha_i&\lambda_i\end{array}\right]$.
             The logarithmic norm inducted
            by the euclidean norm of this matrix, is calculated easily, and equals to $\lambda_i$. We present explicit estimates that were
            used in actual computations in \cite{Supplement}.

        \end{rem}


\section{Example theorems proved by using presented method}
\label{sec:tabelka}
In Table~\ref{table} and Table~\ref{table2} we present data of several theorems that we managed to prove by using
the presented method.

To obtain the results presented in Table~\ref{table} and Table~\ref{table2} we kept the forcing constant
(it was the same as in Theorem~\ref{thm:main2}),
and we were varying the parameter $\nu$.
The radius of the energy absorbing ball $E_0$ \eqref{eq:absorbingBallRadius} was different for each case.

\begin{table}[h!]
\begin{equation*}
        \begin{array}{|c|c|c|c|c|c|c|c|c|c|c|c|}
        \hline
        \mathbf{id} & \mathbf{\nu} & \mathbf{\int_0^{2\pi}{u_0(x)\,dx}}    & \mathbf{E_0}   &   \mathbf{m} & \mathbf{L^+}& \textbf{1.} & \textbf{2.} & \textbf{3.} & \textbf{4.}\\\hline\hline
        1(\text{Thm.~\ref{thm:main2}}) & \paperExampleNu   &     \paperExampleAzero     &    \paperExampleEZero   &  \paperExampleM    &\paperExampleL  &\paperExampleExecutionTime    & \checkmark  & \checkmark  & \checkmark \\\hline
        2 & 1.9   &   0.6\pi   &    1.352   &   8   &  7.38339e-05  &  365.9  & \checkmark  & \checkmark  & \checkmark \\\hline
        3 & 1.85  &    0       &    1.42607 &   8   &  4.64234e-05  &  933.55  & \checkmark & \checkmark  & \checkmark \\\hline
        4 & 1 & 0 & 4.88072 & 14 & 0.0243972 & -^1 & \checkmark & \checkmark  & \\\hline
        5 & 0.97 & 0 & 5.08197 & 16 & 0.969176 & -^1 & \checkmark & \checkmark  & \\\hline
        6 & 0.85 & 0 & 6.75532 & 12 & 3.72357e+24 & -^1 & \checkmark & & \\\hline
        \end{array}
    \end{equation*}
    \caption{Example results obtained}
    {$^1$ - we do not provide the total execution time, as we could not perform numerical integration in time in those cases.}
    \label{table}
\end{table}

\begin{table}[h!]
\begin{equation*}
        \begin{array}{|c|c|c|c|c|c|c|c|c|c|c|}
        \hline
        \mathbf{id} & \mathbf{\nu} & \mathbf{\int_0^{2\pi}{u_0(x)\,dx}}    & \mathbf{E_0}   &   \mathbf{m} & \mathbf{L^+}&     \textbf{1.}     & \textbf{2.} & \textbf{3.} & \textbf{4.}\\\hline\hline
        1 & 8    & 200\pi & 0.0762613 & 6  & 1.59482e-22 & 806.79 & \checkmark & \checkmark & \checkmark \\\hline
        2 & 0.85 & 40\pi  & 6.75532   & 12 & 0.00709885  & -^1 & \checkmark & \checkmark &  \\\hline
        \end{array}
    \end{equation*}
    \caption{Example results obtained with large $\alpha=\int_0^{2\pi}u_0(x)\,dx$}
    {$^1$ - we do not provide the total execution time, as we could not perform numerical integration in time in those cases.}
    \label{table2}
\end{table}

The meaning of the labels in Table~\ref{table} and
Table~\ref{table2} is as follows, {\bf 1.} is the total execution
time in seconds, {\bf 2.} if the existence of a trapping isolating
segment was established, {\bf 3.} if the periodic solution is
locally attracting, {\bf 4.} if the  periodic solution is
attracting globally, $L^+$ is the upper bound for the Lipschitz constant
of $\mappingSymbol$ -- the time shift by $t_p$. The order of the Taylor method was $6$, time
step length was $0.005$ in all cases.

In some cases, namely in the proofs denoted id 4, 5 in
Table~\ref{table}, and id 2 in Table~\ref{table2} the numerical
integration forward in time of the absorbing set was not
performed, as the calculated absorbing set was too large, all
attempts to integrate it using our algorithm resulted in blow ups
of interval enclosures after a short time. Therefore the step 3 of
Algorithm from Section~\ref{sec:algMainAlgorithm} was not
verified, still, it should be possible also in those cases to
perform successfully the step 3 of Algorithm from
Section~\ref{sec:algMainAlgorithm} by, for instance, applying to
the absorbing set some interval set splitting techniques, and then
integrate separately each small piece. In the proof denoted id 6
in Table~\ref{table} the obtained upper bound for the Lipschitz
constant was $>1$, thus we established just the existence of an
orbit within the trapping isolating segment, without resolving the
question whether this orbit is attracting.

\section{Conclusion}

A method of proving the existence of globally attracting periodic
solutions for a class of dissipative PDEs has been presented. A
detailed case study of the viscous Burgers equation with a
nonautonomous forcing function has been provided. All the rigorous
numerics computer software used is available on-line \cite{Software}.

There are several paths for the future development of the
presented method we will pursue. First, we will investigate the
possibility of obtaining a theoretical result of existence of
attracting orbits, with exponential rate of convergence, for
\eqref{eq:vBEq} with periodic boundary conditions for any forcing,
which is a continuous and bounded function of time.
We will address this topic in our forthcoming papers.

We would like to conclude with brief note about our forthcoming results \cite{CZ}.
We established existence of globally attracting solutions asymptotically for 
large $\alpha$ ($u_0$ integral). Additionally, in this case, we obtained 
bounds for the attracting solution of the form $\alpha + O(\frac{\|f\|}{\alpha})$.
In this work we also considered the whole class of smooth forcings, i.e. 
we dropped the assumption of the finite number of nonzero forcing modes.





\providecommand{\digitsPrecision}{6}
\appendix
\section{Numerical data from proof of Theorem \ref{thm:main2}}
\label{sec:numData}
In this appendix we present the numerical data obtained in the algorithm from Section~\ref{sec:algMainAlgorithm} proving
Theorem \ref{thm:main2}.

Program was programmed in C++ language. Program was executed on Linux 32-bit Intel Core i5-2430M CPU @ 2.40 GHz x 4 machine,
compiled with GCC compiler version 4.7.3, and with the following compiler flags (-O0 -frounding-math -ffloat-store).
\begin{itemize}
  \item The Galerkin projection dimension, $m=8$.
  \item The autonomous forcing $l^2$ energy, $E(\{f_k\})=1.31$.
  \item The nonautonomous forcing $l^2$ energy, $E(\{\widetilde{f}_k\})=1.31$.
  \item Upper bound for the total forcing $l^2$ energy,
    \begin{equation*}
      \max_{\substack{g\in\left\{[f_k]\right\}\\t\in[0,t_p]}}{E(\{g_k+\widetilde{f}_k(t)\})}=4.88072.
    \end{equation*}
  \item The absorbing ball radius $E_0=\paperExampleEZero$.
  \item The Lipschitz constant, $L\leq \paperExampleL$.
  \item The absorbing set, $V\oplus\Theta=$
    \begin{equation*}
    \begin{array}{|c|c|c|}\hline\mathbf{k} & \mathbf{\re{a_k}} & \mathbf{\im{a_k}}\\\hline\hline
1 & 0.0259487+[-0.357457,0.357457] & 0.153052+[-0.238867,0.238867]\\
2 & 0.0449597+[-0.150289,0.150289] & -0.0425997+[-0.143525,0.143525]\\
3 & -0.0178005+[-7.0412,7.0412]10^{-2} & 0.0192868+[-7.07427,7.07427]10^{-2}\\
4 & -1.75301\cdot 10^{-3}+[-2.24348,2.24348]10^{-2} & 1.45325\cdot 10^{-3}+[-2.23378,2.23378]10^{-2}\\
5 & 8.24378\cdot 10^{-4}+[-9.60379,9.60379]10^{-3} & 1.45682\cdot 10^{-4}+[-9.60024,9.60024]10^{-3}\\
6 & -3.07739\cdot 10^{-5}+[-4.3149,4.3149]10^{-3} & -1.0924\cdot 10^{-5}+[-4.30936,4.30936]10^{-3}\\
7 & -4.5983\cdot 10^{-5}+[-2.01048,2.01048]10^{-3} & -1.49242\cdot 10^{-5}+[-2.01036,2.01036]10^{-3}\\
8 & -1.96036\cdot 10^{-6}+[-9.87006,9.87006]10^{-4} & 7.01463\cdot 10^{-6}+[-9.86863,9.86863]10^{-4}\\
9 & 2.81191\cdot 10^{-6}+[-6.19244,6.19244]10^{-4} & 4.15566\cdot 10^{-6}+[-6.19241,6.19241]10^{-4}\\
10 & -2.16324\cdot 10^{-7}+[-4.00182,4.00182]10^{-4} & -1.24755\cdot 10^{-6}+[-3.99861,3.99861]10^{-4}\\
11 & -1.10967\cdot 10^{-7}+[-2.84633,2.84633]10^{-4} & -1.0639\cdot 10^{-7}+[-2.84634,2.84634]10^{-4}\\
12 & -5.47851\cdot 10^{-8}+[-2.18552,2.18552]10^{-4} & 1.37878\cdot 10^{-7}+[-2.18508,2.18508]10^{-4}\\
13 & 1.29273\cdot 10^{-8}+[-1.77235,1.77235]10^{-4} & -1.69395\cdot 10^{-8}+[-1.77237,1.77237]10^{-4}\\
14 & 9.27549\cdot 10^{-9}+[-1.49524,1.49524]10^{-4} & -1.9964\cdot 10^{-9}+[-1.4952,1.4952]10^{-4}\\
15 & -4.70792\cdot 10^{-9}+[-1.29781,1.29781]10^{-4} & -2.81077\cdot 10^{-10}+[-1.29781,1.29781]10^{-4}\\
16 & 4.35414\cdot 10^{-10}+[-1.15017,1.15017]10^{-4} & -2.87508\cdot 10^{-10}+[-1.15016,1.15016]10^{-4}\\
17 & 4.76886\cdot 10^{-10}+[-1.03539,1.03539]10^{-4} & 2.82898\cdot 10^{-10}+[-1.03539,1.03539]10^{-4}\\
18 & -2.05696\cdot 10^{-10}+[-9.43268,9.43268]10^{-5} & -8.29394\cdot 10^{-11}+[-9.43267,9.43267]10^{-5}\\
19 & -4.29042\cdot 10^{-12}+[-8.67209,8.67209]10^{-5} & -1.99325\cdot 10^{-11}+[-8.67208,8.67208]10^{-5}\\
20 & 1.68368\cdot 10^{-11}+[-8.02637,8.02637]10^{-5} & 2.14338\cdot 10^{-11}+[-8.02637,8.02637]10^{-5}\\
21 & -1.01567\cdot 10^{-11}+[-7.46033,7.46033]10^{-5} & -1.19795\cdot 10^{-11}+[-7.46033,7.46033]10^{-5}\\
22 & 1.67491\cdot 10^{-12}+[-6.93165,6.93165]10^{-5} & 4.91813\cdot 10^{-12}+[-6.93165,6.93165]10^{-5}\\
23 & 2.14631\cdot 10^{-13}+[-6.41116,6.41116]10^{-5} & -1.35591\cdot 10^{-12}+[-6.41116,6.41116]10^{-5}\\
\geq 24 & \multicolumn{2}{|c|}{|a_k|\leq 11362.2/k^{5},\ |a_{24}|=  -2.41124\cdot 10^{-14}+[-5.86026,5.86026]10^{-5}} \\\hline\end{array}
    \end{equation*}
  \item The trapping region, $\trappingRegion=$
  \begin{equation*}
\begin{array}{|c|c|c|}\hline\mathbf{k} & \mathbf{\re{a_k}} & \mathbf{\im{a_k}}\\\hline\hline
1 & 0.0862439+[-1,1]10^{-4} & 0.11949+[-1,1]10^{-4}\\
2 & 0.0376504+[-1,1]10^{-4} & -0.0422059+[-1,1]10^{-4}\\
3 & -0.0192285+[-1,1]10^{-4} & 0.021161+[-1,1]10^{-4}\\
4 & -2.03091\cdot 10^{-4}+[-1,1]10^{-4} & 5.43697\cdot 10^{-4}+[-1,1]10^{-4}\\
5 & 1.56511\cdot 10^{-4}+[-1,1]10^{-4} & -1.51924\cdot 10^{-5}+[-1,1]10^{-4}\\
6 & -2.91075\cdot 10^{-5}+[-1,1]10^{-4} & 1.7875\cdot 10^{-6}+[-1,1]10^{-4}\\
7 & -1.7401\cdot 10^{-6}+[-1,1]10^{-4} & 4.24224\cdot 10^{-7}+[-1,1]10^{-4}\\
8 & 2.88967\cdot 10^{-7}+[-1,1]10^{-4} & 2.38365\cdot 10^{-7}+[-1,1]10^{-4}\\
9 & 0+[-6.24295,6.24295]10^{-5} & 0+[-6.24295,6.24295]10^{-5}\\
10 & 0+[-4.096,4.096]10^{-5} & 0+[-4.096,4.096]10^{-5}\\
11 & 0+[-2.79762,2.79762]10^{-5} & 0+[-2.79762,2.79762]10^{-5}\\
12 & 0+[-1.97531,1.97531]10^{-5} & 0+[-1.97531,1.97531]10^{-5}\\
13 & 0+[-1.43412,1.43412]10^{-5} & 0+[-1.43412,1.43412]10^{-5}\\
14 & 0+[-1.06622,1.06622]10^{-5} & 0+[-1.06622,1.06622]10^{-5}\\
15 & 0+[-8.09086,8.09086]10^{-6} & 0+[-8.09086,8.09086]10^{-6}\\
16 & 0+[-6.25,6.25]10^{-6} & 0+[-6.25,6.25]10^{-6}\\
17 & 0+[-4.90416,4.90416]10^{-6} & 0+[-4.90416,4.90416]10^{-6}\\
18 & 0+[-3.90184,3.90184]10^{-6} & 0+[-3.90184,3.90184]10^{-6}\\
19 & 0+[-3.14301,3.14301]10^{-6} & 0+[-3.14301,3.14301]10^{-6}\\
20 & 0+[-2.56,2.56]10^{-6} & 0+[-2.56,2.56]10^{-6}\\
21 & 0+[-2.10612,2.10612]10^{-6} & 0+[-2.10612,2.10612]10^{-6}\\
22 & 0+[-1.74851,1.74851]10^{-6} & 0+[-1.74851,1.74851]10^{-6}\\
23 & 0+[-1.46369,1.46369]10^{-6} & 0+[-1.46369,1.46369]10^{-6}\\
\geq 24 & \multicolumn{2}{|c|}{|a_k|\leq 0.4096/k^{4},\ |a_{24}|=  0+[-1.23457,1.23457]10^{-6}} \\\hline\end{array}
    \end{equation*}

  \item Approximate eigenvalues of $D\,P_mG\left(\MID\left(\left[\varphi\left([t_0, t_1], 0,\trappingRegion\right)\right]\right)\right)$ matrix \eqref{eq:diagonalForm},
    $\spect\left(D\,P_mG\left(\MID\left(\left[\varphi\left([t_0, t_1], 0,\trappingRegion\right)\right]\right)\right)\right)=$
    \begin{gather}
    \left\{
    -127.956+i4.00124,-127.956-i4.00124,-98.0028+i3.49945,-98.0028-i3.49945,\nonumber\right.\\
    -2.00765+i0.496188,-2.00765-i0.496188,-8.00731+i0.999835,-8.00731-i0.999835,\nonumber\\
    -18.0068+i1.49997,-18.0068-i1.49997,-32.0066+i1.99999,-32.0066-i1.99999,\nonumber\\
    \left.-72.006+i2.99993,-72.006-i2.99993,-50.0065+i2.49999,-50.0065-i2.49999\label{eq:approxEigenv}\right\}.
    \end{gather}
\end{itemize}

The absorbing set is apparently larger than the trapping region, it was
necessary for the proof to integrate it rigorously forward in time. The Taylor method used in
the $C^0$ Lohner nonautonomous integrator was of order $6$ with time step $0.005$. Total execution
time was $\paperExampleExecutionTime$ seconds.

Here we presented data limited to \digitsPrecision, more detailed numerical data with higher precision is available on-line at \cite{Software}.

\end{document}